\documentclass[10pt]{amsart}

\usepackage{epic,eepic}
\usepackage{latexsym}
\usepackage{amssymb}
\usepackage{graphicx}

\usepackage{latexsym}
\usepackage{amssymb}

\newcommand{\excise}[1]{}

\numberwithin{equation}{section}
\newtheorem{thm}{Theorem}[section]
\newtheorem{lemma}[thm]{Lemma}
\newtheorem{claim}[thm]{Claim}

\newtheorem{prop}[thm]{Proposition}

\newtheorem{theorem}{Theorem}

\newtheorem{Example}[thm]{Example}
\newtheorem{Alg}[thm]{Algorithm}
\newtheorem{Defn}[thm]{Definition}

\newenvironment{remark}{\begin{trivlist}\item {\bf
        Remark.\,}}{\mbox{}\hfill$\square$\end{trivlist}}

\newenvironment{remarks}{\begin{trivlist}\item {\bf
        Remarks.\,}}{\mbox{}\hfill$\square$\end{trivlist}}

    {\begin{eqnarray}}
    {\end{eqnarray}$\!\!$}

        {\begin{list}
                {\noindent\makebox[0mm][r]{\arabic{enumi}.}}
                {\leftmargin=5.5ex \usecounter{enumi}}
        }
        {\end{list}}

\newenvironment{romanlist}%
        {\begin{list}
                {\noindent\makebox[0mm][r]{(\roman{enumi})}}
                {\leftmargin=5.5ex \usecounter{enumi}}
        }
        {\end{list}}

\newenvironment{romanlist'}%
        {\begin{list}
                {\noindent\makebox[0mm][r]{(\roman{enumi}$'$)}}
                {\leftmargin=5.5ex \usecounter{enumi}}
        }
        {\end{list}}

\def\<{\langle}
\def\>{\rangle}
\def\0{{\mathbf 0}}
\def\1{{\mathbf 1}}

\def\CC{{\mathbb C}}

\def\OO{{\mathcal O}}
\def\PP{{\mathbb P}}

\def\mm{{\mathfrak m}}
\def\pp{{\mathbf p}}
\def\qq{{\mathbf q}}

\def\KK{{\mathcal K}}

\def\hh{{\mathfrak h}}

\def\ii{{\bf i}}
\def\gg{{\mathfrak g}}

\def\id{{\rm id}}

\def\proj{{\rm Proj}}
\def\Gr{{\mathcal G}}
\def\Hom{{\rm Hom}}

\def\gln{{G_{\!}L_n}}

\def\cprime{$'$}

\newcommand{\lie}{\mbox{Lie}}
\newcommand{\bfC}{\mathbb{C}}
\newcommand{\cplx}{\mathbb{C}}

\def\II{\mathbb I}
\def\JJ{\mathbb J}
\def\highest{{highest }}
\def\lowest{{lowest }}

\newcommand{\sll}{\rm S_{\!}L}
\newcommand{\gll}{\rm G_{\!}L}

\newcommand{\sln}{\mathit{SL}_n}

\newcommand{\spfour}{\mathit{Sp}_4}
\newcommand{\lgr}{\mathcal{G}}
\newcommand{\lgrgl}{\lgr_{\gll}}
\newcommand{\lgrsl}{\lgr_{\sll}}
\newcommand{\base}{\delta}
\newcommand{\kost}{\mathbb{K}}
\newcommand{\coweights}{\Lambda}
\newcommand{\cowgl}{\Lambda_{\gll}}
\newcommand{\cowsl}{\Lambda_{\sll}}

\newcommand{\dimz}{\dim_0}
\newcommand{\unip}{N}
\newcommand{\momm}{\Phi}

\newcommand{\wcomp}{M}
\newcommand{\comp}{M^\dag}
\newcommand{\scomp}{M^\ddag}

\newcommand{\polytinv}{\mathbb{C}[t,t^{-1}]}
\newcommand{\polyt}{\mathbb{C}[t]}

\newcommand{\integers}{\mathbb{Z}}

\newcommand{\cdim}{\mbox{dim}}

\newcommand{\eqz}{{\equiv_{{\text{\tiny $0$}}}}}
\newcommand{\len}{{\rm len}}
\newcommand{\highweight}{\lambda}
\newcommand{\kst}{\pp}

\begin{document}

\title[Mirkovi\'c-Vilonen cycles and polytopes in type A]%
    {Mirkovi\'c-Vilonen cycles and polytopes in type A}
\author{Jared Anderson}
\address{University of Pittsburgh\\Pittsburgh, PA\\USA}
\email{jared@math.pitt.edu}

\author{Mikhail Kogan}
\thanks{MK was supported by the NSF Postdoctoral Fellowship}
\address{Northeastern University\\Boston, MA\\USA}
\email{misha@research.neu.edu}


\begin{abstract}
\noindent We study, in type A, the algebraic cycles (MV-cycles)
discovered by I.~Mirkovi\'c and K.~Vilonen~\cite{MV}.  In
particular, we partition the loop Grassmannian into smooth pieces
such that the MV-cycles are their closures.  We explicitly
describe the points in each piece using the lattice model of the
loop Grassmannian in type A.  The partition is invariant under the
action of the coweights and, up to this action, the pieces are
parametrized by the Kostant parameter set.  This description of MV-cycles
allows us to prove the main result of the paper: the computation of the
moment map images of MV-cycles (MV-polytopes) by identifying the vertices
of each polytope.

\vspace{3mm}

\noindent
(Mathematics subject classification number: 14L35)

\end{abstract}

\maketitle

\section{Introduction}
\label{sec:intro}

Let $G$ be a connected complex algebraic reductive group. The loop
Grassmannian $\Gr$ for $G$ is the quotient $G(\KK)/G(\OO)$ where
$\OO=\CC[[t]]$ is the ring of formal power series and $\KK=\CC((t))$ is
its field of fractions, the ring of formal Laurent series. The same set is
obtained using polynomials instead of power series:
$\Gr=G(\CC[t,t^{-1}])/G(\CC[t])$. $\lgr$ may be realized as an increasing
union of finite-dimensional complex algebraic varieties~\cite{L}; more
formally, it has the structure of an ind-scheme whose set of geometric
points is $G(\KK)/G(\OO)$. In type~A, when $G=\sln$ or $G=\gln$, there is
a simple, well-known lattice model for the loop Grassmannian, described in
Section~\ref{sec:lattice-model}. It allows one to study $\Gr$ using
finite-dimensional linear algebra, as we will do throughout.

Our main results are concerned with certain finite-dimensional
algebraic subvarieties of $\Gr$, which we call MV-cycles.  These
were discovered by Mirkovi\'c and Vilonen~\cite{MV} and provide a
canonical basis of algebraic cycles for the intersection homology
of the closures of the strata of $\Gr$, for the natural
stratification by $G(\OO)$-orbits.

To define MV-cycles, we first need some more notation. Choose a
maximal torus $T$ of $G$ and let $\coweights=\Hom (\CC^*,T)$ be
the coweight lattice. A coweight $\lambda$, viewed as a map from
$\CC^*$ to $T$, is given by a trigonometric series, which, after
substituting $t$ for $e^{i\theta}$, produces an element in
$G(\KK)$, and hence in $\Gr$, which we denote by $\underline
\lambda$.  These are the $T$-fixed points in $\Gr$. A coweight
$\lambda$, treated as an element of $G(\KK)$, defines a
\emph{shift} $\sigma_\lambda$ of $\Gr$ by left multiplication.
This is an isomorphism of $\Gr$, taking any fixed point
$\underline \alpha$ to the fixed
point~$\underline{\alpha+\lambda}$.

Denote by $N$ and $N^-$ the unipotent radicals of opposite Borel subgroups
of $G$ intersecting along~$T$. Each $N(\KK)$-orbit on $\Gr$ contains a
unique $T$-fixed point $\underline{\alpha}$; similarly for $N^-(\KK)$. For
$\alpha,\beta\in \coweights$, set $S_{\alpha,\beta} = N(\KK)
\underline{\alpha} \cap N^-(\KK) \underline{\beta}$; this is nonempty if
and only if $\alpha-\beta$ is a linear combination of positive coroots
with nonnegative integer coefficients.  We define the \emph{MV-cycles}
with \highest coweight $\alpha$ and \lowest coweight $\beta$ to be the
irreducible components of the closure of $S_{\alpha,\beta}$. As shown
in~\cite{A2}, the results of~\cite{MV} imply that this definition of
MV-cycles is equivalent (up to shifts) to the one in~\cite{MV}. 
Except for this equivalence of definitions, this paper is 
independent of the results of~\cite{MV}. 
In particular, we provide, in type A, an independent proof of the
important purity property of MV-cycles established in~\cite{MV}, by
showing that  $S_{\alpha,\beta}$ is a complex algebraic
variety of pure complex dimension equal to the height of $\alpha-\beta$.


As shown in \cite{AP}, the left action of the maximal compact
torus $T_K\subset T$ on $\Gr$ can be viewed as a Hamiltonian
action whose moment map values lie in ${\rm Lie}(T_K)$. Given an
MV-cycle $M$, we denote its image under this moment map by $P_M$
and call it an \mbox{\emph{MV-polytope}}. See~\cite{A2} for many
examples of MV-polytopes in low rank groups.

Our aim is to provide a detailed, explicit, and useful description
and parametrization of \mbox{MV-cycles} and MV-polytopes in type A. This
is not straightforward because of the use of irreducible
components in the definition. I.~Mirkovi\'c and M.~Finkelberg,
using work of G.~Lusztig, have a different way of finding a dense
open subset of each cycle, which works for any type~\cite{M}.
The results of Gaussent and Littelmann~\cite{GL}, which we 
learned of when this paper was almost completed, also provide 
a different explicit description of MV-cycles, which works for any type. 
It would be very interesting to compare these approaches, 
particularly with a view to calculating MV-polytopes in all types.

There are many reasons for wanting to understand MV-cycles as
explicitly as possible; let us mention three of them.  First,
these varieties are intrinsically interesting since they provide a
canonical basis for representations of the (Langlands dual) group.
Second, as shown in~\cite{A2}, MV-polytopes allow for
combinatorial calculations in representation theory: weight
multiplicities and tensor product multiplicities equal the number
of MV-polytopes fitting inside a certain region. Third, as
conjectured in~\cite{A1}, there is a natural Hopf algebra
structure on the vector space spanned by (equivalence classes
under shifts of) MV-cycles, isomorphic to the algebra of functions
on the group $N$. Moreover there seems to be a canonical set of
generators and relations, which were calculated in the cases of
$\spfour$ and $\sln$, $n\leq 5$, by looking at MV-polytopes. As
observed by B.~Sturmfels and A.~Zelevinsky, these provide examples
of cluster algebras in the sense of Fomin and
Zelevinsky~\cite{FZ2}. We hope that the explicit understanding of
MV-cycles here will eventually allow for the construction of this
cluster algebra structure in type A.

Using the lattice model for the loop  Grassmannian in type A,
namely when $G=\sln$ or $G=\gln$, we provide, in
Theorem~\ref{thm:main}, a decomposition of $\Gr$ into smooth
pieces whose closures are the MV-cycles. The partition is
invariant under the action of the coweights and, up to this
action, the pieces are parametrized by the Kostant parameter set.
The description of the pieces
is very explicit and in the spirit of the definition of Bruhat
cells in the Grassmannian: just as every Bruhat cell contains
those vector spaces for which the dimensions of intersections with
a fixed flag is determined by a partition, every piece of our
decomposition contains those lattices for which dimensions
(understood in the proper sense) of intersections with a fixed
collection of infinite-dimensional vector spaces is determined by
a Kostant partition.

Using this decomposition, we prove our main result: the computation 
of MV-polytopes by explicitly
identifying their vertices in Theorem~\ref{thm:main2}.  The key is
a combinatorial algorithm presented in Section~\ref{algorithm},
which, given an MV-polytope, constructs certain lower-dimensional
MV-polytopes that are faces of it.
Repeated application of the
algorithm yields the vertices.
While the algorithm is combinatorial in nature, we don't know how
to prove some of its properties without using the geometry of the
loop Grassmannian, our decomposition in particular.

The paper is largely self-contained, 
using only linear algebra,
combinatorics, and the most basic complex algebraic geometry.  
One exception is that we take as given, in
Section~\ref{sec:decompositions}, the existence and basic
properties of the moment map for the torus action.

We begin, in Section~\ref{sec:main-theorem}, by describing the
lattice model for $\Gr$ in type A and use it to state our main
results: Theorem~\ref{thm:main} and Theorem~\ref{thm:main2}. In
Section~\ref{sec:combinatorics}, we discuss certain combinatorial
properties of the abovementioned algorithm.
Section~\ref{sec:bases} is the heart of the paper and uses the
lattice model to construct open dense sets of points in each
cycle.  Included here is an explicit identification
(Proposition~\ref{prop:M-is-smooth}) of each piece of our
decomposition with an open Zariski subset of a product of
projective spaces. In Section~\ref{sec:decompositions}, which is
the only one containing results for all types, we discuss the
moment map and the moment images of three well-known coarse
decompositions of the loop Grassmannian by certain group actions.
Theorem~\ref{thm:cones} describes the general shape of moment
images of subvarieties of~$\Gr$. In
Section~\ref{sec:moment-images}, we discuss the relation of our
decomposition of the loop Grassmannian to the finer decomposition
into torus orbit types. We calculate moment images of torus orbits
and, in Theorem~\ref{thm:strong-moment}, identify all torus orbits
inside each MV-cycle for which the moment map image of the orbit's
closure is the whole MV-polytope. We end this section by proving
our main results.

\smallskip

\noindent{\bf Acknowledgments.} Both authors thank R.~MacPherson,
I.~Mirkovi\'c, M.~Vybornov and A.~Zelevinsky for helpful discussions and
suggestions. MK is especially indebted to A. Zelevinsky for introducing
him to the subject of cluster algebras and motivating this work by
pointing to the connection between MV-cycles and cluster algebras.

\section{Main results} \label{sec:main-theorem}

In this section we give precise statements of our main results,
which are in type A.  We partition the loop Grassmannian into
smooth pieces such that the MV-cycles are the closures of these
pieces.  This decomposition is invariant under shifts by
coweights, and, up to shifts, the pieces are parametrized by the
Kostant parameter set.  We give a combinatorial algorithm for
constructing the vertices of the MV-polytope corresponding to an
element of the Kostant parameter set.

\subsection{Lattice model}
\label{sec:lattice-model}

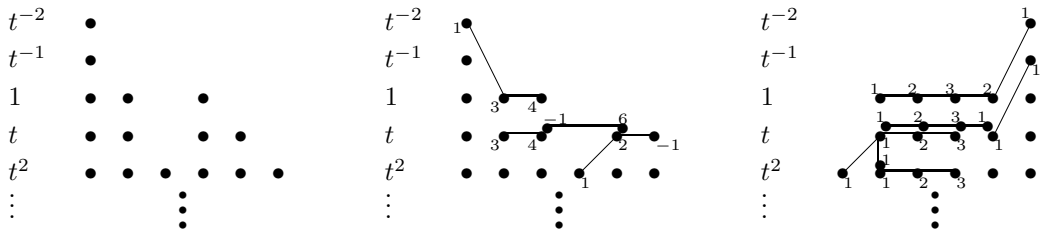
\begin{figure}[ht]
\begin{picture}(400,73)

\put(0,10){ \setlength{\unitlength}{.5mm}

\def\smalldot{{$\bullet$}}
\def\tinydot{{\tiny$\bullet$}}

\put(0,0){$t^2$} \put(0,10){$t$} \put(0,20){$1$}
\put(0,30){$t^{-1}$} \put(0,40){$t^{-2}$} \put(0,-10){$\vdots$}

\put(20,0){\smalldot} \put(30,0){\smalldot} \put(40,0){\smalldot}
\put(50,0){\smalldot} \put(60,0){\smalldot} \put(70,0){\smalldot}

\put(20,10){\smalldot} \put(30,10){\smalldot}
\put(50,10){\smalldot} \put(60,10){\smalldot}

\put(20,20){\smalldot} \put(30,20){\smalldot}
\put(50,20){\smalldot}

\put(20,30){\smalldot} \put(20,40){\smalldot}

\put(45,-5){\tinydot} \put(45,-9){\tinydot} \put(45,-13){\tinydot}


\put(100,0){$t^2$} \put(100,10){$t$} \put(100,20){$1$}
\put(100,30){$t^{-1}$} \put(100,40){$t^{-2}$}
\put(100,-10){$\vdots$}

\put(120,0){\smalldot} \put(130,0){\smalldot}
\put(140,0){\smalldot} \put(150,0){\smalldot}
\put(160,0){\smalldot} \put(170,0){\smalldot}

\put(120,10){\smalldot} \put(130,10){\smalldot}
\put(140,10){\smalldot} \put(160,10){\smalldot}
\put(170,10){\smalldot}

\put(120,20){\smalldot} \put(130,20){\smalldot}
\put(140,20){\smalldot}

\put(120,30){\smalldot} \put(120,40){\smalldot}


\put(132,22.5){\line(1,0){10}} \put(132,12.5){\line(1,0){10}}
\put(132,22){\line(-1,2){10}} \put(143,14.5){\line(1,0){20}}
\put(152,2){\line(1,1){10}} \put(172,12){\line(-1,0){10}}


\put(141.5,12){\smalldot} \put(161.5,12){\smalldot}


{\tiny

\put(118,39){$1$} \put(128,18){$3$} \put(138,18){$4$}

\put(128,8){$3$} \put(138,8){$4$}

\put(152,-2){$1$} \put(162,8){$2$} \put(172,8){$-1$}

\put(142,15){$-1$} \put(162,15){$6$}}

\put(145,-5){\tinydot} \put(145,-9){\tinydot}
\put(145,-13){\tinydot}


\put(200,0){$t^2$} \put(200,10){$t$} \put(200,20){$1$}
\put(200,30){$t^{-1}$} \put(200,40){$t^{-2}$}
\put(200,-10){$\vdots$}

\put(220,0){\smalldot} \put(230,0){\smalldot}
\put(240,0){\smalldot} \put(250,0){\smalldot}
\put(260,0){\smalldot} \put(270,0){\smalldot}

 \put(230,10){\smalldot}
\put(240,10){\smalldot} \put(250,10){\smalldot}
\put(260,10){\smalldot} \put(270,10){\smalldot}

\put(230,20){\smalldot} \put(240,20){\smalldot}
\put(250,20){\smalldot} \put(260,20){\smalldot}
\put(270,20){\smalldot}

\put(270,30){\smalldot} \put(270,40){\smalldot}


\put(232,22.5){\line(1,0){30}} \put(232,12.5){\line(1,0){20}}
\put(232,14.5){\line(1,0){28}} \put(232,2.5){\line(1,0){20}}

\put(222,2){\line(1,1){10}}

\put(262,12){\line(1,2){10}} \put(262,22){\line(1,2){10}}

\put(231.2,2){\line(0,1){10}}


\put(231.5,12.5){\smalldot} \put(241.5,12.5){\smalldot}
\put(251.5,12.5){\smalldot} \put(258.5,12.5){\smalldot}

\put(230,2){\smalldot}


{\tiny

\put(232,-2){$1$}\put(242,-2){$2$}\put(252,-2){$3$}

\put(232,8.5){$1$} \put(242,8){$2$} \put(252,8){$3$}
\put(232,4){$1$} \put(222,-2){$1$}

\put(230.5,15.5){$1$}\put(240.5,15.5){$2$}\put(250.5,15.5){$3$}
\put(257.5,15.5){$1$}

\put(229,23){$1$} \put(239,23){$2$} \put(249,23){$3$}
\put(259,23){$2$} \put(269,43){$1$}

\put(262,8){$1$} \put(272,28){$1$}

\put(245,-5){\tinydot} \put(245,-9){\tinydot}
\put(245,-13){\tinydot} }

}

\end{picture}
\caption{Examples of lattices.} \label{fig:lattice-examples}
\end{figure}

Following Lusztig~\cite{L},
we will use the lattice model of the loop Grassmannian~$\lgr$ in
type~A.  Points in $\lgr$ are subspaces of a certain
infinite-dimensional complex vector space, satisfying some extra
conditions.  The vector space $X$ is defined by specifying a
basis: if $e_1, e_2, \ldots, e_n$ denotes the standard basis
for~$\cplx^n$, then $X$ is the span of $t^je_i$ for $1\leq i\leq
n$, $-\infty<j<\infty$, where we regard this as a symbol with two
indices.  We usually picture these basis vectors in an array
consisting of $n$ columns, infinite in both directions. Let $t$ be
the invertible linear operator on $X$ that sends each~$t^je_i$ to
$t^{j+1}e_i$.  Then $\lgrgl$ consists of those subspaces $Y$ of
$X$ such that
\begin{enumerate} \item $Y$ is closed under the action of $t$.
\item $t^NX_0\subseteq Y\subseteq t^{-N}X_0$ for some $N$,
  where $X_0$ is the span of those $t^je_i$ with $j\geq0$.
\end{enumerate}
We call such $Y$ {\em lattices}. Three examples of lattices are
illustrated by the pictures in Figure~\ref{fig:lattice-examples}.
Each dot, or set of dots connected by line segments, represents a
vector; these, together with all the dots $t^je_i$, $j\geq 3$
below the pictures are a $\CC$-basis for the lattice. The first
picture represents the torus fixed point associated, by our
convention, to the coweight $(2,0,-2,0,-1,-2)$. The second picture
represents the lattice that is generated as a $\CC[t]$-module by
the vectors $t^{-2}e_1+3e_2+4e_3$, $3te_2+4te_3$, $-te_3+6te_5$,
$t^2e_4+2te_5-te_6$, $t^{2}e_5$ and $t^{2}e_6$. (Of course, this
set of generators is not unique; for example $t^2e_5$ can be
replaced by~$t^2e_3$.)
The third picture represents the lattice that is generated as a
$\CC[t]$-module by the vectors $e_2+2e_3+3e_4+2e_5+t^{-2}e_6$,
$t^2 e_1+t^2e_2 +te_2+2te_3+3te_4$, $te_5+t^{-1}e_6$, $e_6$,
$t^3e_2$ and~$t^3e_3$.

The space $\lgrgl$ is the loop Grassmannian for $\gln$: one checks
that~$\gln(\polytinv)$ acts transitively on lattices and that the
stabilizer of $X_0$ is $\gln(\polyt)$. The action is the obvious
one suggested by the notation: $\gln (\bfC)$ acts on
$e_1,\ldots,e_n$ by the standard representation, and letting this
commute with $t$ provides an action of $\gln(\polytinv)$ on $X$,
which induces an action on lattices.

The {\em relative dimension} of a lattice $Y$ with respect to
$X_0$ is defined to be $\cdim(Y/Y\cap X_0)-\cdim(X_0/Y\cap X_0)$.
The lattices in Figure~\ref{fig:lattice-examples} have relative
dimensions $-3$, $-6$ and $-9$. The connected components of
$\lgrgl$ are parametrized by relative dimension and are isomorphic
to each other by appropriate shifts $\sigma_\lambda$. The
component consisting of the lattices with relative dimension~$0$
is the loop Grassmannian $\lgrsl \subset \lgrgl$ for $\sln$: one
checks that~$\sln(\polytinv)$ acts transitively on such lattices
and that the stabilizer of $X_0$ is $\sln(\polyt)$.

Let $\cowgl=\integers^n$ and
$\cowsl=\{(\lambda_1,\dots,\lambda_n)\in \cowgl \mid \sum
\lambda_i=0 \}$ denote the coweight lattices of~$\gln$ and~$\sln$.
For any lattice $Y$, define $\base(Y)$ to be the $n$-tuple of
integers $(\base_1(Y),\dots,\base_n(Y))$ where~$\base_i(Y)$ is the
largest $j$ such that $t^{-j}e_i \in Y$.  Let $Y_0$ denote the
vector space spanned by all $t^{-j}e_i\in Y$ and let $\dimz Y
=\dim(Y/Y_0)$. Notice that for $Y\in \lgrgl$ we have $\base(Y) \in
\cowgl$; for $Y\in \lgrsl$ we have $\base(Y) \in \cowsl$ if and
only if $Y=Y_0$. Conversely, for each coweight $\lambda$, the
lattice $\underline{\lambda}$
is spanned by all~$t^{-j}e_i$ with~$j\leq \lambda_i$. We choose
the maximal torus consisting of diagonal matrices, and opposite
Borels consisting of upper and lower triangular matrices. Then the
dominant coweights are those~$\lambda$ for which $i_1<i_2$ implies
$\lambda_{i_1}\geq \lambda_{i_2}$; the coweight is called
\emph{strictly} dominant if $i_1<i_2$ implies~$\lambda_{i_1}>
\lambda_{i_2}$.

\subsection{Kostant parameter set}
\label{sec:kostant-parameter}

Now we define the parameter set.  Recall that the Dynkin diagram
for $\sln$ is $n-1$ dots in a row, one for each simple root
(connected by line segments, which we will not draw). We denote a
positive root by a loop around a sequence of consecutive dots in
the Dynkin diagram, and we call the number of dots it encloses the
length of the loop.  The simple roots are the loops of length $1$.
The other positive roots are loops of length $\geq2$; each is the
sum of the simple roots corresponding to the dots it encloses. (Of
course the word {\em loop} here has nothing to do with its use in
{\em loop Grassmannian}.) A {\em Kostant picture} is a picture of
the Dynkin diagram together with a finite number of such loops. We
draw the loops so that if the dots contained in one loop are a
proper subset of the dots contained in another, then the one is
encircled by the other; if two loops contain precisely the same
dots, we still draw one encircled by the other. In this way, the
loops in a Kostant picture are partially ordered by encirclement.
We write $L\subset L'$ if loop~$L'$ encircles loop~$L$, and
$L\subseteq L'$ if either $L=L'$ or $L'$ encircles~$L$.

\begin{figure}[ht]
\begin{center}
\begin{picture}(380,35)
\setlength{\unitlength}{.5mm}

\def\smalldot{{$\bullet$}}

\put(0,10){\smalldot} \put(15,10){\smalldot}
\put(30,10){\smalldot} \put(45,10){\smalldot}
\put(60,10){\smalldot}

\put(100,10){\smalldot} \put(115,10){\smalldot}
\put(160,10){\smalldot} \put(130,10){\smalldot}
\put(145,10){\smalldot}

\put(200,10){\smalldot} \put(215,10){\smalldot}
\put(260,10){\smalldot} \put(230,10){\smalldot}
\put(245,10){\smalldot}


\thicklines


\put(30,10){\put(1.5,1.6){\ellipse{6}{6}}}
\put(22.5,10){\put(1.5,1.6){\ellipse{25}{10}}}
\put(37.5,10){\put(1.5,1.6){\ellipse{25}{10}}}
\put(52.5,10){\put(1.5,1.6){\ellipse{28}{10}}}
\put(52.5,10){\put(1.5,1.6){\ellipse{25}{7}}}
\put(52.5,10){\put(1.5,1.6){\ellipse{31}{13}}}
\put(30,10){\put(1.5,1.6){\ellipse{44}{15}}}

\put(115,10){\put(1.5,1.6){\ellipse{6}{6}}}
\put(152.5,10){\put(1.5,1.6){\ellipse{28}{9}}}
\put(137.5,10){\put(1.5,1.6){\ellipse{28}{9}}}
\put(107.5,10){\put(1.5,1.6){\ellipse{28}{11}}}


\put(260,10){\put(1.5,1.6){\ellipse{6}{6}}}
\put(222.5,10){\put(1.5,1.6){\ellipse{23}{8}}}
\put(237.5,10){\put(1.5,1.6){\ellipse{63}{18}}}
\put(230,10){\put(1.5,1.6){\ellipse{42}{13}}}
\put(215,10){\put(1.5,1.6){\ellipse{42}{13}}}

\end{picture}
\end{center}
\caption{Examples of Kostant pictures.}
\label{fig:Kostant-pictures}
\end{figure}
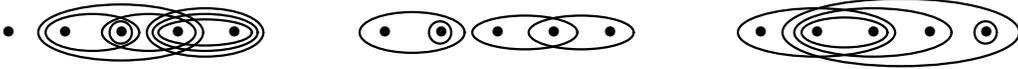

Define the \emph{length} $\len(\pp)$ of a Kostant picture $\pp$ to
be the sum of the lengths of the loops in~$\pp$, and~$|\pp|$ to be
the number of loops in~$\pp$. The {\em Kostant parameter set}
$\kost$ is the collection of all Kostant pictures. Examples of
Kostant pictures for $n=6$ are in
Figure~\ref{fig:Kostant-pictures}. If $\alpha_{ij}$ is the root
that is the sum of simple roots $\alpha_i+\dots+\alpha_{j-1}$,
then these three pictures represent the following Kostant
partitions:
$\alpha_{3}+\alpha_{24}+\alpha_{35}+\alpha_{25}+3\alpha_{46}$,
$\alpha_{2}+\alpha_{13}+\alpha_{35}+\alpha_{46}$ and
$\alpha_{5}+\alpha_{24}+\alpha_{14}+\alpha_{25}+\alpha_{26}$.

It will be useful to imagine each of the $n-1$ dots of the Dynkin
diagram as lying on the boundary between two of the $n$ columns of
basis vectors for $X$.  Then, associated to each loop $L$, will be
the vector subspace $V_L$ spanned by the columns the loop passes
through. To be precise, for each $i=1,2,...,n$, let $V_i$ be the
span of the basis vectors $t^je_i$ in the $i^{\rm th}$ column. If we
number the dots of the Dynkin diagram $1,2,...,n-1$ and if a loop
$L$ contains dots $\ell,\ell+1,...,r-1$ then $V_L=V_\ell\oplus
V_{\ell+1}\oplus ... \oplus V_{r}$.  We say that the loop $L$ has
its {\em left end} at column~$\ell$, its {\em right end} at
column~$r$, and {\em passes through} columns $\ell,\ell+1,...,r$.

\subsection{Parametrization of MV-cycles}
\label{sec:MV-parametrization} We partition $\lgrgl$ into pieces
parametrized by $\kost \times \cowgl$ by defining a map which
sends $Y\in \lgrgl$ to $(\kst(Y),\highweight(Y))\in \kost \times
\cowgl$ as follows.
Let~$\kst(Y)$ be the Kostant picture in which the number of loops
encircling the same dots as a loop $L$ is
\begin{align*}
n_L=&\dimz(Y\cap V_L) - \dimz(Y\cap V_{L_+}) - \dimz(Y\cap
V_{L_-})+\dimz(Y\cap V_{L_{\pm}}) \\
=&\dim\Big((Y\cap V_L)\big/\big((Y\cap V_{L_+})+(Y\cap
V_{L_-})\big)\Big).
\end{align*}
In the above expressions,
$L_+$,~$L_-$ and~$L_{\pm}$ are the loops obtained from $L$
by removing the leftmost,
rightmost, or both of these dots respectively. (If $L$ has length
$1$ then $V_{L_+}, V_{L_-}$ should be taken as the zero vector
space; if $L$ has length $1$ or $2$, then~$V_{L_\pm}$ should be
taken as the zero vector space.) Also, we are viewing $Y\cap
V_L$ as a lattice in $V_L$, as we shall throughout the paper.
We let
$\highweight(Y)=(\lambda_1,\dots,\lambda_n)$ be the coweight such
that $\lambda_i=\base_i(Y)+l_i$ where $l_i$ is the number of loops
in~$\kst(Y)$ whose left end is at column $i$. Then the loop
Grassmannian $\lgrgl$ decomposes into
pieces~$\wcomp(\pp,\lambda)$, each of which contains those
lattices $Y$ with $(\kst(Y),\highweight(Y))=(\pp,\lambda)$. Notice
that both $\underline{\highweight(Y)}$ and~$Y$ have the same
relative dimension~$|\pp|$ with respect to $Y_0$; therefore if
$Y\in \lgrsl$ then $\highweight(Y)\in \cowsl$.
So the parametrization restricts from $\cowgl$ to $\cowsl$: we
have $Y\in \lgrsl \mapsto (\kst(Y),\highweight(Y))\in \kost \times
\cowsl$, and if $\lambda \in \cowsl$ then $\wcomp(\pp,\lambda)
\subset \lgrsl$.

A related decomposition of $\lgrgl$ is the one into torus orbit
types: let $\Gr^T$ be the set of torus fixed points, and, for a
subset $S$ of $\Gr^T$, let $M_S$ be
the set of lattices $Y$ such that $S$ is the set of fixed points
contained in the closure of the torus orbit through $Y$.
The loop Grassmannian is the disjoint union of all
nonempty~$M_S$'s.
Our first theorem identifies MV-cycles as the closures of the
pieces $\wcomp(\pp,\lambda)$ and states that the decomposition of
the loop Grassmannian into the $M_S$'s is a refinement of the
decomposition into the $\wcomp(\pp,\lambda)$'s.

\begin{theorem}
\label{thm:main} Every $\wcomp(\pp,\lambda)$ is a smooth algebraic
variety of dimension $\len(\pp)$ and is a union of a finite number
of~$M_S$'s.
The MV-cycles with highest coweight $\lambda$ are
$\{\overline{\wcomp(\pp,\lambda)} \mid \pp \in \kost
\}$.
\end{theorem}

\begin{remarks}

\noindent (1) In Proposition~\ref{prop:M-is-smooth}, we will
canonically identify each $\wcomp(\pp,\lambda)$ with an open
Zariski subset of a~product of $|\pp|$ projective spaces.

\smallskip

\noindent (2) The usual Grassmannian of $k$-planes in $\cplx^n$
also has a decomposition into torus orbit types, defined just as
above; this is discussed in~\cite{GGMS}. It is still unknown, even
for the Grassmannian, which~$M_S$ are nonempty. In the case of the
loop Grassmannian, using Theorems~\ref{thm:main2}
and~\ref{thm:cones}, it will be easy to identify certain sets~$S$
for which~$M_S$ is nonempty.

The $M_S$ decomposition is invariant under the action of the
affine Weyl group, but the $\wcomp(\pp,\lambda)$ decomposition is
of course not; it is invariant only under the affine part, namely
translations. But it is easy to construct from it a finer,
Weyl-invariant decomposition: for every element~$w$ of the
(ordinary) Weyl group~$W$ there exists a decomposition of the loop
Grassmannian into pieces $M^w(\pp,\lambda)$, each of which is
$\wcomp(\pp,\lambda)$ acted upon by $w$.  Intersecting these $|W|$
decompositions---so that each piece has the form $\bigcap_{w\in W}
M^w(\pp,\lambda)$---is a Weyl-invariant decomposition of the loop
Grassmannian. In the remark after Lemma~\ref{lem:moment}, we will
show that this decomposition is identical with the
$M_S$~decomposition.

\smallskip

\noindent(3) Using Theorem~\ref{thm:main} together with the definition of
MV-cycles, it is easy to identify all lattices in an orbit
$N(\KK)\underline \lambda$. The resulting description of the orbit
was originally obtained in \cite{Ngo}. 
\end{remarks}

\subsection{Vertices of MV-polytopes}
\label{algorithm} Fix $(\pp,\lambda)\in\kost \times \cowgl$;
by Theorem~\ref{thm:main} there
corresponds an MV-cycle $M=\overline{\wcomp(\pp,\lambda)}$ and,
applying the moment map, an MV-polytope $P_M=P(\pp,\lambda)$. We
will describe $P_M$ by constructing a map from the Weyl group $W$
onto its vertices; $W$ is the group of permutations of
$\{1,2,...,n\}$ and the vertices are coweights.

Fix $w\in W$. We will give an $n$-step construction of the
corresponding vertex $\nu=\nu(w)$. In terms of the lattice
$\underline{\nu}$, the idea is essentially simple: starting with
column $w(1)$ we put in it as many basis vectors as possible,
according to the length of the longest chain of nested loops in
$\pp$ passing through that column; then, after removing these
loops from $\pp$ in a certain sense, we put as many basis vectors
as possible in column $w(2)$; and so on up to column $w(n)$.  But
it will take a little work to state this precisely.

\begin{figure}[ht]
\begin{center}
\def\smalldot{{\put(-80,-80){$\bullet$}}}
\def\dot{{\put(-100,-100){\Large$\bullet$}}}

\setlength{\unitlength}{0.00041667in}
\begingroup\makeatletter\ifx\SetFigFont\undefined%
\gdef\SetFigFont#1#2#3#4#5{%
  \reset@font\fontsize{#1}{#2pt}%
  \fontfamily{#3}\fontseries{#4}\fontshape{#5}%
  \selectfont}%
\fi\endgroup%
{\renewcommand{\dashlinestretch}{30}
\begin{picture}(5959,8021)(0,-10)

\thicklines

\put(9014,3576){\ellipse{480}{480}}

\put(8377,3576){\smalldot} \put(9013,3576){\smalldot}

\put(8632,3574){\ellipse{1438}{940}}
\put(8371,3576){\ellipse{658}{658}}
\put(8372,3576){\ellipse{480}{480}}
\put(8378,3576){\ellipse{300}{300}}
\put(9013,3573){\ellipse{300}{300}}

\path(4888,7522)(4888,7222) \path(4888,7522)(4888,7222)
\thicklines \texture{44000000 aaaaaa aa000000 8a888a 88000000
aaaaaa aa000000 888888
    88000000 aaaaaa aa000000 8a8a8a 8a000000 aaaaaa aa000000 888888
    88000000 aaaaaa aa000000 8a888a 88000000 aaaaaa aa000000 888888
    88000000 aaaaaa aa000000 8a8a8a 8a000000 aaaaaa aa000000 888888 }
\shade\path(2570,5867)(4564,4707)(4564,2387)
    (3567,1807)(577,3547)(577,4707)(2570,5867)
\path(2570,5867)(4564,4707)(4564,2387)
    (3567,1807)(577,3547)(577,4707)(2570,5867)
\path(-2000,4000)(-2000,3000)(-1333,2500)
\path(-2100,3800)(-2000,4000)(-1900,3800)
\path(-1500,2500)(-1333,2500)(-1413,2700) \put(-2500,
3500){$\alpha_1$} \put(-1900, 2500){$\alpha_2$}
\thinlines \path(1573,5287)(1573,2967) \path(3567,5287)(3567,1807)
\put(2571,5834){\dot} \put(582,3550){\dot} \put(592,4721){\dot}
\put(4561,4698){\dot} \texture{80222222 22555555 55808080 80555555
55222222 22555555 55880888 8555555
    55222222 22555555 55808080 80555555 55222222 22555555 55080808 8555555
    55222222 22555555 55808080 80555555 55222222 22555555 55880888 8555555
    55222222 22555555 55808080 80555555 55222222 22555555 55080808 8555555 }
\put(4883,3254){\smalldot} \put(4883,3254){\smalldot}
\put(4883,2790){\smalldot} \put(4883,2340){\smalldot}
\put(4883,1890){\smalldot} \put(4883,1440){\smalldot}
\put(5330,1440){\smalldot} \put(5780,1440){\smalldot}
\put(5780,1890){\smalldot} \put(5780,2340){\smalldot}
\put(5708,945){\makebox(0,0)[lb]{\smash{{{\SetFigFont{14}{16.8}{\rmdefault}{\mddefault}{\updefault}2}}}}}
\put(5258,945){\makebox(0,0)[lb]{\smash{{{\SetFigFont{14}{16.8}{\rmdefault}{\mddefault}{\updefault}3}}}}}
\put(4808,945){\makebox(0,0)[lb]{\smash{{{\SetFigFont{14}{16.8}{\rmdefault}{\mddefault}{\updefault}1}}}}}
\put(83,1710){\smalldot} \put(530,1710){\smalldot}
\put(980,1710){\smalldot} \put(980,2160){\smalldot}
\put(980,2610){\smalldot} \put(530,2160){\smalldot}
\put(530,2610){\smalldot} \put(983,3060){\smalldot}
\put(530,3060){\smalldot}
\put(24,1200){\makebox(0,0)[lb]{\smash{{{\SetFigFont{14}{16.8}{\rmdefault}{\mddefault}{\updefault}3}}}}}
\put(474,1200){\makebox(0,0)[lb]{\smash{{{\SetFigFont{14}{16.8}{\rmdefault}{\mddefault}{\updefault}2}}}}}
\put(924,1200){\makebox(0,0)[lb]{\smash{{{\SetFigFont{14}{16.8}{\rmdefault}{\mddefault}{\updefault}1}}}}}
\put(83,5580){\smalldot} \put(530,5580){\smalldot}
\put(980,5580){\smalldot} \put(980,6030){\smalldot}
\put(980,6480){\smalldot} \put(530,6030){\smalldot}
\put(530,6480){\smalldot} \put(530,6930){\smalldot}
\put(530,7380){\smalldot}
\put(23,5070){\makebox(0,0)[lb]{\smash{{{\SetFigFont{14}{16.8}{\rmdefault}{\mddefault}{\updefault}3}}}}}
\put(473,5070){\makebox(0,0)[lb]{\smash{{{\SetFigFont{14}{16.8}{\rmdefault}{\mddefault}{\updefault}1}}}}}
\put(923,5070){\makebox(0,0)[lb]{\smash{{{\SetFigFont{14}{16.8}{\rmdefault}{\mddefault}{\updefault}2}}}}}
\put(3995,6434){\makebox(0,0)[lb]{\smash{{{\SetFigFont{12}{14.4}{\rmdefault}{\mddefault}{\updefault}t}}}}}
\put(4115,6584){\makebox(0,0)[lb]{\smash{{{\SetFigFont{8}{9.6}{\rmdefault}{\mddefault}{\updefault}-1}}}}}
\put(3999,6914){\makebox(0,0)[lb]{\smash{{{\SetFigFont{12}{14.4}{\rmdefault}{\mddefault}{\updefault}t}}}}}
\put(4119,7064){\makebox(0,0)[lb]{\smash{{{\SetFigFont{8}{9.6}{\rmdefault}{\mddefault}{\updefault}-2}}}}}
\put(3992,5084){\makebox(0,0)[lb]{\smash{{{\SetFigFont{12}{14.4}{\rmdefault}{\mddefault}{\updefault}t}}}}}
\put(4112,5234){\makebox(0,0)[lb]{\smash{{{\SetFigFont{8}{9.6}{\rmdefault}{\mddefault}{\updefault}2}}}}}
\put(4003,5556){\makebox(0,0)[lb]{\smash{{{\SetFigFont{12}{14.4}{\rmdefault}{\mddefault}{\updefault}t}}}}}
\put(3973,6006){\makebox(0,0)[lb]{\smash{{{\SetFigFont{12}{14.4}{\rmdefault}{\mddefault}{\updefault}1}}}}}
\path(4338,7006)(4638,7006) \path(4338,7006)(4638,7006)
\path(4338,6106)(4638,6106) \path(4338,6106)(4638,6106)
\path(4338,5206)(4638,5206) \path(4338,5206)(4638,5206)
\path(4338,5656)(4638,5656) \path(4338,5656)(4638,5656)
\path(4338,6556)(4638,6556) \path(4338,6556)(4638,6556)
\path(5788,7522)(5788,7222) \path(5788,7522)(5788,7222)
\path(5338,7522)(5338,7222) \path(5338,7522)(5338,7222)
\put(4561,2384){\dot}
\put(4787,7664){\makebox(0,0)[lb]{\smash{{{\SetFigFont{12}{14.4}{\rmdefault}{\mddefault}{\updefault}e}}}}}
\put(4937,7589){\makebox(0,0)[lb]{\smash{{{\SetFigFont{8}{9.6}{\rmdefault}{\mddefault}{\updefault}1}}}}}
\put(5683,7657){\makebox(0,0)[lb]{\smash{{{\SetFigFont{12}{14.4}{\rmdefault}{\mddefault}{\updefault}e}}}}}
\put(5833,7582){\makebox(0,0)[lb]{\smash{{{\SetFigFont{8}{9.6}{\rmdefault}{\mddefault}{\updefault}3}}}}}
\put(5252,7664){\makebox(0,0)[lb]{\smash{{{\SetFigFont{12}{14.4}{\rmdefault}{\mddefault}{\updefault}e}}}}}
\put(5402,7589){\makebox(0,0)[lb]{\smash{{{\SetFigFont{8}{9.6}{\rmdefault}{\mddefault}{\updefault}2}}}}}
\put(4883,7021){\smalldot} \put(4883,6557){\smalldot}
\put(4883,6107){\smalldot} \put(4883,5657){\smalldot}
\put(4883,5207){\smalldot} \put(5330,5207){\smalldot}
\put(5780,5207){\smalldot} \put(5330,6107){\smalldot}
\put(5330,5657){\smalldot}
\put(4813,4693){\makebox(0,0)[lb]{\smash{{{\SetFigFont{14}{16.8}{\rmdefault}{\mddefault}{\updefault}1}}}}}
\put(5263,4693){\makebox(0,0)[lb]{\smash{{{\SetFigFont{14}{16.8}{\rmdefault}{\mddefault}{\updefault}2}}}}}
\put(5713,4693){\makebox(0,0)[lb]{\smash{{{\SetFigFont{14}{16.8}{\rmdefault}{\mddefault}{\updefault}3}}}}}
\put(3143,1879){\smalldot} \put(3143,1429){\smalldot}
\put(3143,979){\smalldot} \put(3143,529){\smalldot}
\put(4040,529){\smalldot} \put(4040,979){\smalldot}
\put(4040,1429){\smalldot} \put(4040,1879){\smalldot}
\put(3590,529){\smalldot}
\put(3069,0){\makebox(0,0)[lb]{\smash{{{\SetFigFont{14}{16.8}{\rmdefault}{\mddefault}{\updefault}2}}}}}
\put(3519,0){\makebox(0,0)[lb]{\smash{{{\SetFigFont{14}{16.8}{\rmdefault}{\mddefault}{\updefault}3}}}}}
\put(3969,0){\makebox(0,0)[lb]{\smash{{{\SetFigFont{14}{16.8}{\rmdefault}{\mddefault}{\updefault}1}}}}}
\put(2288,6509){\smalldot} \put(2735,6509){\smalldot}
\put(3185,6509){\smalldot} \put(2735,6959){\smalldot}
\put(2735,7409){\smalldot} \put(2285,6959){\smalldot}
\put(2285,7409){\smalldot} \put(2735,7859){\smalldot}
\put(2735,8323){\smalldot}
\put(2214,5999){\makebox(0,0)[lb]{\smash{{{\SetFigFont{14}{16.8}{\rmdefault}{\mddefault}{\updefault}2}}}}}
\put(2664,5999){\makebox(0,0)[lb]{\smash{{{\SetFigFont{14}{16.8}{\rmdefault}{\mddefault}{\updefault}1}}}}}
\put(3114,5999){\makebox(0,0)[lb]{\smash{{{\SetFigFont{14}{16.8}{\rmdefault}{\mddefault}{\updefault}3}}}}}
\put(3571,1821){\dot} \path(1573,5287)(4564,3547)
\path(4564,3547)(2570,2387) \path(577,3547)(3567,5287)
\path(2570,5867)(2570,2387) \path(4564,4707)(1573,2967)
\path(577,4707)(4564,2387)

\end{picture}
}
\end{center}
\caption {An example of an MV-polytope} \label{fig:mv}
\end{figure}
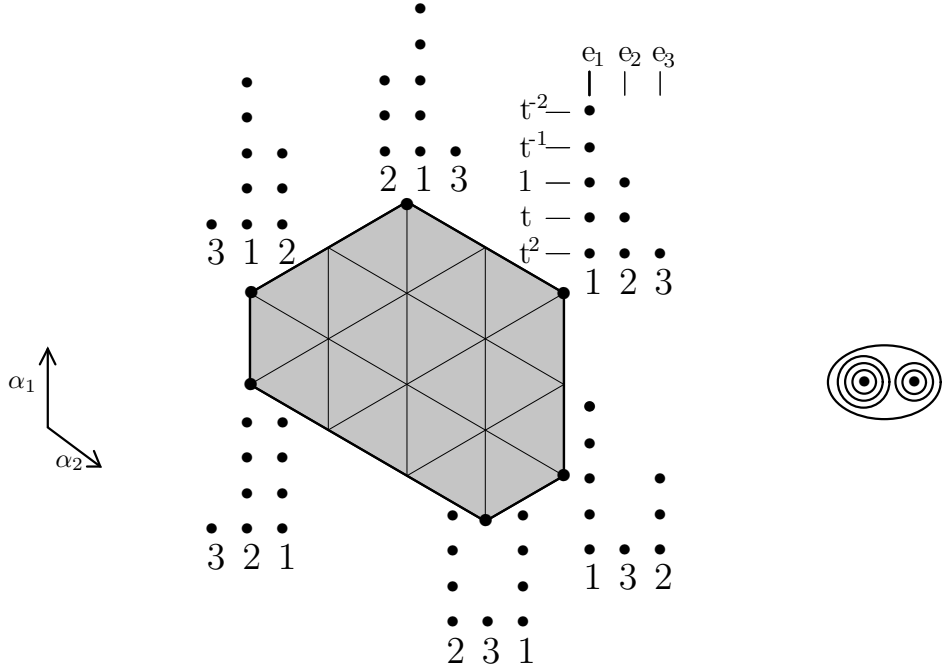

An example is drawn in Figure~\ref{fig:mv}: the Kostant picture
$\pp$ is drawn to the right of the corresponding MV-polytope; the
directions of the two simple roots are shown on the left; the
vertices $\nu(w)$ are labelled by the pictures of the lattices
$\underline{\nu(w)}$, with Weyl group element~$w$ written below
(an $i$ written below column $j$ means $w(i)=j$); the top right
vertex is the highest coweight $\lambda$.

To define $\nu(w)$, we first construct, for any Kostant picture
$\pp$ in type $A_{r-1}$ and for any column~$i$ $(1\leq i\leq r)$,
a Kostant picture $\hat \pp_i$ in type $A_{r-2}$, called the
\emph{collapse} of $\pp$ along column $i$.
Let $\qq$ be the Kostant picture consisting of all loops in $\pp$
that pass through column $i$, and rank these loops by level of
encirclement: level~1 consists of all those loops in $\qq$ that
don't encircle any loops in $\qq$; level~2 consists of all those
loops in $\qq$ that don't encircle any loops in $\qq$ except those
in level~1; level~3 consists of all those loops in $\qq$ that
don't encircle any loops in $\qq$ except those in levels~1 and~2;
etc.

Fix an arbitrary level. Notice that the loops in this level do not
encircle each other and thus are naturally ordered from left to
right, say $L_1, L_2,..., L_k$. Form a new sequence with one fewer
loop, $M_1, M_2,..., M_{k-1}$ as follows: each~$M_i$ is the
{\em{join}} of loops $L_i$ and $L_{i+1}$, meaning that the left
end of~$M_i$ is the left end of~$L_i$ and the right end of~$M_i$
is the right end of $L_{i+1}$. (Of course if $k=1$ then this new
sequence is empty.)

Do this for each level, and let $\qq^\prime$ be the Kostant
picture consisting of all the new loops formed for all levels. Let
$\pp^\prime$ be the Kostant picture obtained by replacing $\qq$ by
$\qq^\prime$; that is, $\pp^\prime$ consists of all loops that are
in $\pp$ but not in $\qq$ together with all loops in $\qq^\prime$.
It is easy to see that no loop of $\pp^\prime$ has its left end or
its right end at column $i$. So if we simply cross out column $i$
and identify the two dots of the Dynkin diagram on either side, we
see that $\pp^\prime$ may be viewed as a Kostant picture in one
rank lower, namely type $A_{r-2}$. (In fact it is the Kostant
picture corresponding to a particular face of the MV-polytope.)
This is the Kostant picture $\hat \pp_i$. Notice that $|\hat
\pp_i|$ equals $|\pp|$ minus the number of levels.

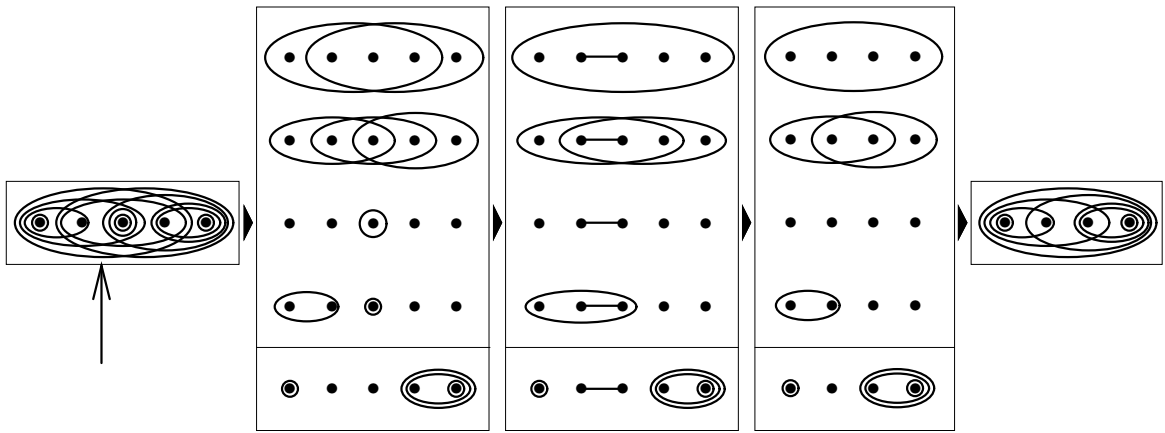
\begin{figure}[ht]
\begin{center}
\def\smalldot2{{\put(-115,-115){$\bullet$}}}

\setlength{\unitlength}{0.00028333in}
\begingroup\makeatletter\ifx\SetFigFont\undefined%
\gdef\SetFigFont#1#2#3#4#5{%
  \reset@font\fontsize{#1}{#2pt}%
  \fontfamily{#3}\fontseries{#4}\fontshape{#5}%
  \selectfont}%
\fi\endgroup%
{\renewcommand{\dashlinestretch}{30}
\begin{picture}(20887,7689)(0,-10)
\put(2112,3762){\smalldot2} \thicklines
\put(2112,3762){\ellipse{300}{300}}
\put(2112,3763){\ellipse{480}{480}} \thinlines
\put(1365,3762){\smalldot2} \put(2862,3762){\smalldot2}
\thicklines \put(3612,3762){\ellipse{300}{300}} \thinlines
\put(3612,3762){\smalldot2} \put(609,3762){\smalldot2} \thicklines
\put(612,3762){\ellipse{300}{300}}
\put(912,3762){\ellipse{1160}{540}}
\put(3292,3762){\ellipse{1160}{540}}
\put(3292,3762){\ellipse{1340}{700}}
\put(1378,3762){\ellipse{2260}{840}}
\put(2122,3762){\ellipse{2260}{840}}
\put(2877,3762){\ellipse{2260}{1000}}
\put(2505,3762){\ellipse{3200}{1240}}
\put(1760,3762){\ellipse{3200}{1240}} \thinlines
\path(12,4512)(4212,4512)(4212,3012)
    (12,3012)(12,4512)
\put(18783,3762){\smalldot2} \put(18036,3762){\smalldot2}
\put(19533,3762){\smalldot2} \thicklines
\put(20283,3762){\ellipse{300}{300}} \thinlines
\put(20283,3762){\smalldot2} \thicklines
\put(19963,3762){\ellipse{1160}{540}}
\put(19963,3762){\ellipse{1340}{700}}
\put(18793,3762){\ellipse{2260}{840}}
\put(19548,3762){\ellipse{2260}{1000}}
\put(19176,3762){\ellipse{3200}{1240}}
\put(18033,3762){\ellipse{300}{300}}
\put(18347,3762){\ellipse{1160}{540}} \thinlines
\put(6627,6742){\smalldot2} \put(5880,6742){\smalldot2}
\put(7377,6742){\smalldot2} \put(8127,6742){\smalldot2}
\put(5124,6742){\smalldot2} \thicklines
\put(6275,6742){\ellipse{3200}{1240}} \thinlines
\put(6627,5242){\smalldot2} \put(5880,5242){\smalldot2}
\put(7377,5242){\smalldot2} \put(8127,5242){\smalldot2}
\put(5124,5242){\smalldot2} \thicklines
\put(5893,5242){\ellipse{2260}{840}}
\put(6637,5242){\ellipse{2260}{840}}
\put(7392,5242){\ellipse{2260}{1000}} \thinlines
\put(6627,3742){\smalldot2} \thicklines
\put(6627,3743){\ellipse{480}{480}} \thinlines
\put(5880,3742){\smalldot2} \put(7377,3742){\smalldot2}
\put(8127,3742){\smalldot2} \put(5124,3742){\smalldot2}
\put(6627,2242){\smalldot2} \thicklines
\put(6627,2242){\ellipse{300}{300}} \thinlines
\put(5880,2242){\smalldot2} \put(7377,2242){\smalldot2}
\put(8127,2242){\smalldot2} \put(5124,2242){\smalldot2}
\thicklines \put(5427,2242){\ellipse{1160}{540}}
\put(7020,6742){\ellipse{3200}{1240}} \thinlines
\put(11130,6742){\smalldot2} \put(10383,6742){\smalldot2}
\put(11880,6742){\smalldot2} \put(12630,6742){\smalldot2}
\put(9627,6742){\smalldot2} \put(11130,5242){\smalldot2}
\put(10383,5242){\smalldot2} \put(11880,5242){\smalldot2}
\put(12630,5242){\smalldot2} \put(9627,5242){\smalldot2}
\put(11130,3742){\smalldot2} \put(10383,3742){\smalldot2}
\put(11880,3742){\smalldot2} \put(12630,3742){\smalldot2}
\put(9627,3742){\smalldot2} \put(11130,2242){\smalldot2}
\put(10383,2242){\smalldot2} \put(11880,2242){\smalldot2}
\put(12630,2242){\smalldot2} \put(9627,2242){\smalldot2}
\thicklines \put(11142,6742){\ellipse{4000}{1240}}
\put(10735,5242){\ellipse{3000}{840}}
\put(11499,5242){\ellipse{3000}{840}}
\put(10384,2242){\ellipse{2000}{580}} \thinlines
\put(14913,6762){\smalldot2} \put(14166,6762){\smalldot2}
\put(15663,6762){\smalldot2} \put(16413,6762){\smalldot2}
\thicklines \put(15306,6762){\ellipse{3200}{1240}} \thinlines
\put(14913,5262){\smalldot2} \put(14166,5262){\smalldot2}
\put(15663,5262){\smalldot2} \put(16413,5262){\smalldot2}
\thicklines \put(14923,5262){\ellipse{2260}{840}}
\put(15678,5262){\ellipse{2260}{1000}} \thinlines
\put(14913,3762){\smalldot2} \put(14166,3762){\smalldot2}
\put(16413,3762){\smalldot2} \put(14913,2262){\smalldot2}
\put(14166,2262){\smalldot2} \put(15663,2262){\smalldot2}
\put(16413,2262){\smalldot2} \thicklines
\put(14477,2262){\ellipse{1160}{540}} \thinlines
\put(15663,3762){\smalldot2} \put(14913,772){\smalldot2}
\put(14166,772){\smalldot2} \put(15663,772){\smalldot2}
\thicklines \put(16413,772){\ellipse{300}{300}} \thinlines
\put(16413,772){\smalldot2} \thicklines
\put(16093,772){\ellipse{1160}{540}}
\put(16093,772){\ellipse{1340}{700}}
\put(14163,772){\ellipse{300}{300}} \thinlines
\put(6627,757){\smalldot2} \put(5880,757){\smalldot2}
\put(7377,757){\smalldot2} \thicklines
\put(8127,757){\ellipse{300}{300}} \thinlines
\put(8127,757){\smalldot2} \put(5124,757){\smalldot2} \thicklines
\put(5127,757){\ellipse{300}{300}}
\put(7807,757){\ellipse{1160}{540}}
\put(7807,757){\ellipse{1340}{700}} \thinlines
\put(11130,757){\smalldot2} \put(10383,757){\smalldot2}
\put(11880,757){\smalldot2} \thicklines
\put(12630,757){\ellipse{300}{300}} \thinlines
\put(12630,757){\smalldot2} \put(9627,757){\smalldot2} \thicklines
\put(9630,757){\ellipse{300}{300}}
\put(12310,757){\ellipse{1160}{540}}
\put(12310,757){\ellipse{1340}{700}} \path(1727,1227)(1727,2982)
\path(1877.000,2382.000)(1727.000,2982.000)(1577.000,2382.000)
\path(10375,762)(11125,762) \path(10375,2262)(11125,2262)
\path(10375,3762)(11125,3762) \path(10375,5262)(11125,5262)
\path(10375,6762)(11125,6762) \thinlines
\path(4525,1512)(8725,1512) \drawline(8725,1512)(8725,1512)
\path(9025,1512)(13225,1512) \path(13525,1512)(17125,1512)
\path(17425,4512)(20875,4512)(20875,3012)
    (17425,3012)(17425,4512)
\path(4525,7662)(8725,7662)(8725,12)
    (4525,12)(4525,7662)
\path(9025,7662)(13225,7662)(13225,12)
    (9025,12)(9025,7662)
\path(13525,7662)(17125,7662)(17125,12)
    (13525,12)(13525,7662)
\blacken\path(4300,4062)(4300,3462)(4450,3762)
    (4300,4062)(4300,4062)
\path(4300,4062)(4300,3462)(4450,3762)
    (4300,4062)(4300,4062)
\blacken\path(13300,4062)(13300,3462)(13450,3762)
    (13300,4062)(13300,4062)
\path(13300,4062)(13300,3462)(13450,3762)
    (13300,4062)(13300,4062)
\blacken\path(17200,4062)(17200,3462)(17350,3762)
    (17200,4062)(17200,4062)
\path(17200,4062)(17200,3462)(17350,3762)
    (17200,4062)(17200,4062)
\blacken\path(8804,4062)(8804,3462)(8954,3762)
    (8804,4062)(8804,4062)
\path(8804,4062)(8804,3462)(8954,3762)
    (8804,4062)(8804,4062)
\end{picture}
}

\end{center}
\caption{Construction of $\hat \pp_i$ out of $\pp$.}
\label{fig:kstep}
\end{figure}

This process is illustrated in Figure~\ref{fig:kstep} for the
Kostant picture $\pp$ at the left of the figure and $i=3$ as
indicated by the arrow; $\hat \pp_3$ is drawn on the right; $\qq$
contains eight of the twelve loops, which have been separated into
four levels in the second column of the figure.

Now we successively collapse the columns in the order determined
by $w$. For an ordered subset $\II=(i_1\dots, i_k)$ of
$\{1,\dots,n\}$ let $\hat \pp_\II$ denote the Kostant picture that
results from collapsing columns $i_1\dots,i_k$ in this order.
(Actually, since after each collapse there is one fewer column,
and the remaining columns are renumbered, we should, to be precise,
say it like this:
Let $h(m)$ be the number of $j<m$ with $i_j<i_m$. Then
$\hat\pp_\II$ is produced by collapsing column~$i_1$, then
column~$i_2-h(2)$, then column~$i_3-h(3)$ and so on.)

Now we can define $\nu=\nu(w)$.  As before, let $l_i$ denote the
number of loops in $\pp$ whose left end is at column $i$. For
$1\leq k\leq n$, let $\II_k=(w(1),\dots,w(k))$, and let
$N_{w(k)}=|\hat\pp_{\II_k}|-|\hat\pp_{\II_{k-1}}|$ denote the
number of loops removed at the $k^{\rm th}$ step during the
collapse along column~$w(k)$. Set
\begin{equation}
\label{eq:def_nu}
\nu_i=\lambda_i-l_i+N_i.
\end{equation}
Note that if $\lambda \in \cowsl$ then $\nu\in \cowsl$ since $\sum
\nu_i-\lambda_i=|\pp|-\sum l_i =0$.

Notice that if $w$ is the identity permutation then
$\nu(w)=\lambda$; this is the highest coweight vertex of the
polytope.

\begin{figure}[ht]
\begin{center}
\setlength{\unitlength}{0.00053333in}

\def\smalldot3{{\put(-60,-60){$\bullet$}}}

\begingroup\makeatletter\ifx\SetFigFont\undefined%
\gdef\SetFigFont#1#2#3#4#5{%
  \reset@font\fontsize{#1}{#2pt}%
  \fontfamily{#3}\fontseries{#4}\fontshape{#5}%
  \selectfont}%
\fi\endgroup%
{\renewcommand{\dashlinestretch}{30}
\begin{picture}(10366,2003)(0,-10)
\path(1369,12)(1369,762)
\path(1444.000,582.000)(1369.000,762.000)(1294.000,582.000)
\path(4849,12)(4849,762)
\path(4924.000,582.000)(4849.000,762.000)(4774.000,582.000)
\path(8212,12)(8212,762)
\path(8287.000,582.000)(8212.000,762.000)(8137.000,582.000)
\path(9562,12)(9562,762)
\path(9637.000,582.000)(9562.000,762.000)(9487.000,582.000)
\path(9967,12)(9967,762)
\path(10042.000,582.000)(9967.000,762.000)(9892.000,582.000)
\put(4548,1379){\smalldot3} \put(3953,1379){\smalldot3}
\put(5145,1379){\smalldot3} \put(5742,1379){\ellipse{238}{238}}
\put(5742,1379){\smalldot3} \put(5487,1379){\ellipse{924}{430}}
\put(5487,1379){\ellipse{1066}{556}}
\put(4556,1379){\ellipse{1798}{668}}
\put(5157,1379){\ellipse{1798}{796}}
\put(4861,1379){\ellipse{2546}{986}}
\put(3951,1379){\ellipse{238}{238}}
\put(4201,1379){\ellipse{924}{430}} \put(1683,1379){\smalldot3}
\put(1683,1379){\ellipse{238}{238}}
\put(1683,1380){\ellipse{382}{382}} \put(1089,1379){\smalldot3}
\put(2280,1379){\smalldot3} \put(2877,1379){\ellipse{238}{238}}
\put(2877,1379){\smalldot3} \put(487,1379){\smalldot3}
\put(489,1379){\ellipse{238}{238}}
\put(728,1379){\ellipse{924}{430}}
\put(2622,1379){\ellipse{924}{430}}
\put(2622,1379){\ellipse{1066}{556}}
\put(1099,1379){\ellipse{1798}{668}}
\put(1691,1379){\ellipse{1798}{668}}
\put(2292,1379){\ellipse{1798}{796}}
\put(1996,1379){\ellipse{2546}{986}}
\put(1403,1379){\ellipse{2546}{986}}
\path(12,1976)(3354,1976)(3354,782)
    (12,782)(12,1976)
\put(7338,1381){\smalldot3} \put(6743,1381){\smalldot3}
\put(7934,1381){\smalldot3} \put(6741,1381){\ellipse{238}{238}}
\put(7934,1381){\ellipse{238}{238}}
\put(7346,1381){\ellipse{1592}{556}}
\put(7336,1381){\ellipse{1718}{668}}
\path(6397,1856)(8307,1856)(8307,902)
    (6397,902)(6397,1856)
\put(8813,1379){\smalldot3} \put(8810,1379){\ellipse{238}{238}}
\put(9407,1379){\smalldot3} \path(8631,1618)(9586,1618)(9586,1140)
    (8631,1140)(8631,1618)
\put(10118,1379){\smalldot3} \put(10116,1379){\ellipse{238}{238}}
\path(9877,1618)(10354,1618)(10354,1140)
    (9877,1140)(9877,1618)
\path(3533,1976)(6159,1976)(6159,782)
    (3533,782)(3533,1976)
\blacken\path(3391,1618)(3391,1140)(3510,1379)(3391,1618)
\path(3391,1618)(3391,1140)(3510,1379)(3391,1618)
\blacken\path(6217,1618)(6217,1140)(6336,1379)(6217,1618)
\path(6217,1618)(6217,1140)(6336,1379)(6217,1618)
\blacken\path(8392,1618)(8392,1140)(8511,1379)(8392,1618)
\path(8392,1618)(8392,1140)(8511,1379)(8392,1618)
\blacken\path(9682,1618)(9682,1140)(9801,1379)(9682,1618)
\path(9682,1618)(9682,1140)(9801,1379)(9682,1618)
\end{picture}
}

\end{center}
\caption{Construction of $\nu(w)$.} \label{fig:kost}
\end{figure}
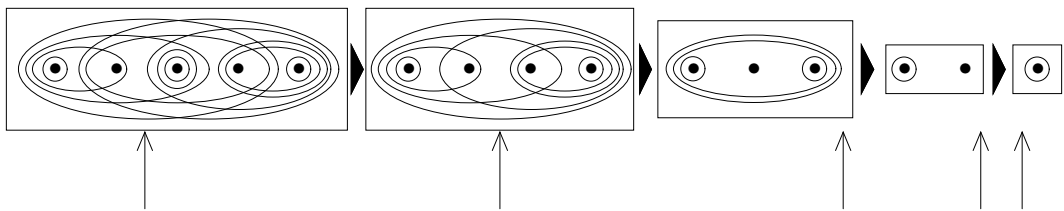

To illustrate how $\nu(w)$ is defined, let $\pp$ be the
leftmost Kostant picture in Figure~\ref{fig:kost}, $\lambda$
the coweight $(2,0,1,0,-1,-2)$, and $w$ the permutation
$(346512)$---that is $w(1)=3$, $w(2)=4$ and so on.
To define $\nu(w)$ we have to perform the collapses shown in
Figure~\ref{fig:kost}.
So $N_{w(1)}=N_3=4$ is the number of
loops removed during the first collapse, $N_{w(2)}=N_4=4$ is the
number of loops removed during the second collapse and so on:
$N_{w(3)}=N_6=3$, $N_{w(4)}=N_5=0$, $N_{w(5)}=N_1=1$, and
$N_{w(6)}=N_2=0$. Hence, by~(\ref{eq:def_nu}), $\nu_1=2-4+1$,
$\nu_2=0-2+0$, and so on, so that $\nu$ is the coweight
$(-1,-2,2,2,-2,1)$.
Note that the first collapse is shown in detail in
Figure~\ref{fig:kstep}.

\begin{theorem}
\label{thm:main2}

The vertices of the MV-polytope $P(\pp,\lambda)$ are $\{ \nu(w)
\mid w \in W \}$.
\end{theorem}

\begin{remark} I.~Mirkovi\'c and M.~Vybornov~\cite{MVy} have found a
decomposition of the loop Grassmannian in type A into the quiver
varieties of Nakajima~\cite{Na1,Na2}.  Like MV-cycles, these
varieties are torus invariant and it would be very interesting to
understand their moment map images. \end{remark}

\section{Combinatorics of collapses}
\label{sec:combinatorics}

This section explains certain combinatorial properties of the
algorithm defined in Section~\ref{algorithm}.
Most proofs are deferred until Section~\ref{sec:bases}
since they rely on geometry.
We end, however, with a purely combinatorial proof of a key
claim used in the proof of Proposition~\ref{strongly-compatible}.

\subsection{Properties of the algorithm} Recall that $\hat \pp_i$
denotes the Kostant picture produced by collapsing~$\pp$ along
column~$i$ and that if $\II$ is the ordered set $(i_1,\dots,i_k)$
then $\hat \pp_\II$ is produced by collapsing along columns
$i_1,\dots, i_k$ in this order.

An unexpected property of this algorithm is that it is commutative:

\begin{theorem} \label{thm:commutativity}
For any Kostant picture~$\pp$ and two columns $i_1,i_2$ we have
$\hat\pp_{(i_1,i_2)}=\hat\pp_{(i_2,i_1)}$.
\end{theorem}

The proof of Theorem~\ref{thm:commutativity} will follow from
certain geometric properties of the decomposition $\wcomp(\pp,\lambda)$,
as explained in the remark preceding
Proposition~\ref{strongly-compatible}.
We would be interested in knowing a combinatorial proof of it.

The next proposition gives some insight into how the algorithm
behaves with multiple collapses and into why
Theorem~\ref{thm:commutativity} is true.
For a Kostant picture $\pp$ and an ordered set $\II$, let $M$ be a
loop of~$\hat\pp_\II$. Define the \emph{ancestry} of $M$ to be the
set of all the loops of $\pp$ used to produce $M$. So if during
the collapse of the last column of $\II$, loop $M$ is produced by
joining two loops $M'$ and $M''$, then the ancestry of~$M$ is the
union of the ancestries of $M'$ and $M''$. The ancestry depends in
general on the order of collapse, as can be seen, for example, by
collapsing the first Kostant picture in
Figure~\ref{fig:Kostant-pictures} along $i=3,4$ in both orders; in
only one of these does the ancestry of the single resulting loop
contain the loop of length~$1$ around the third dot.

Let
$L_1,\dots,L_s$ be the loops from the ancestry of $M$ that do not
encircle any other loop from the ancestry. We say that $M$ is the
\emph{join} of $L_1,\dots, L_s$.  Since the loops $L_1,\dots,L_s$
do not encircle each other, we can assume they are ordered from
left to right: the left and right ends of $L_m$ are to the left of
the corresponding ends of~$L_{m+1}$.

\begin{prop}
\label{prop:s-1}
Suppose a loop $M$ is the join of loops $L_1,\dots, L_s$.
Then the number of columns in the set~$\II$ that are
passed through by at least one of these loops
is exactly $s-1$.
\end{prop}

\begin{remark}
Theorem~\ref{thm:commutativity} and Proposition~\ref{prop:s-1} suggest
that there should be a non-recursive definition of $\hat\pp_\II$,
that describes the simultaneous collapse $\hat\pp_I$ along the
columns $I$; here $I$ is the underlying unordered set of $\II$.
Each loop $M$ of $\hat\pp_I$ should should be a join of loops
$L_1,\dots,L_s$ from $\pp$ not encircling each other,
passing through exactly $s-1$ columns of $I$,
and satisfying additional conditions that we only know
how to state in the case $|I|=1$.
\end{remark}

Proposition~\ref{prop:s-1} follows from
Claim~\ref{claim2}, proved combinatorially below,
and the fact, proved geometrically during the
proof of Proposition~\ref{strongly-compatible},
that we necessarily have $|J|+1=s$ in this claim.

\begin{claim}
\label{claim2}
Suppose a loop $M$ is the join of loops $L_1,\dots, L_s$.
Let $J=(j_1<\dots<j_{|J|}) \subseteq I$ be the set of all columns
passed through by at least one of these loops.
Then $|J|+1\leq s$. Moreover, if $|J|+1=s$, then loop~$L_m$ passes
through columns $j_{m-1}$ and $j_{m}$.
\end{claim}

\subsection{Proof of Claim~\ref{claim2}} We say that a loop $L$ is a \emph{top} loop of a
Kostant picture if it is not encircled by any other loop of this
picture. Without loss of generality we can assume that the top
loops of $\pp$ are the loops $L_1\dots, L_s$. Indeed, removing any
loop of $\pp$ that is not encircled by one of the $L_j$ does not
change the ancestry of $M$. After this, it is clear that we can
also assume that $I=J$; so $J=(j_1<\dots<j_k)$ is a reordering of
the set $\II=(i_1,\dots,i_k)$.

Set $\II_m=(i_1,\dots,i_m)$ for $1\leq m\leq k$; so $\II=\II_k$.
We will prove that the number of top loops of~$\hat\pp_{\II_m}$ is
strictly greater than the number of top loops of
$\hat\pp_{\II_{m+1}}$. Since the Kostant picture $\hat\pp_\II$ has
only one top loop, $M$, this will imply the first part of the
claim.

By induction, it is enough to prove that $\hat
\pp_{\II_1}=\hat\pp_{i_1}$ has fewer top loops than $\pp$. First
consider those loops of $\pp$ that do not pass through column
$i_1$.  Whether or not any particular one of them is a top loop
clearly does not change after the collapse along column $i_1$;
indeed, if it is a top loop, it cannot be encircled by the join of
two loops in the collapse since it does not pass through
column~$i_1$; and if it is not a top loop it cannot be uncovered
by the vanishing of a loop in the collapse since loops in the
ancestry---top loops in particular---do not vanish.

Now consider those loops of $\pp$ passing through column $i_1$.
Recall that the algorithm defining $\hat\pp_{i_1}$ decomposes them
into levels and then joins loops in each level.
Every loop on the highest level is obviously a top loop of~$\pp$.
The algorithm joins loops of the highest level producing top loops
of $\hat\pp_{i_1}$, and the number of loops decreases by one.

Now consider an arbitrary level consisting of loops
$R_1,\dots,R_t$, so that loops $R_m$ and $R_{m+1}$ join to produce
loop~$S_m$ of $\hat\pp_{i_1}$. To finish the proof of the first
part of the claim, we will show that the number of top loops among
$S_1,\dots, S_{t-1}$ is not greater than the number of top loops
among loops $R_1,\dots,R_t$. This  follows immediately from the
following two statements:
\begin{enumerate}
\item{If $R_m$ and $R_{m+1}$ are not top loops of $\pp$ then $S_m$ is not
a top loop of $\hat\pp_{i_1}$.}
\item{If for some $\ell<r-1$,
loops $R_{\ell}$ and $R_r$ are not top loops of~$\pp$, but all the
loops $R_{\ell+1},\dots, R_{r-1}$ are top loops of~$\pp$, then
loops $S_\ell,\dots,S_{r-1}$ are not top loops of~$\hat
\pp_{i_1}$.}
\end{enumerate}

To prove (1) consider rightmost top loop $L_j$ that
encircles~$R_m$ and the leftmost top loop~$L_{j'}$ that encircles
$R_{m+1}$. During the collapse of column~$i_1$ every top loop
passing through column~$i_1$ (in particular loops~$L_j$
and~$L_{j'}$) has to be joined with another loop, since otherwise
this top loop will not be in the ancestry of~$M$. If~$L_j$
encircles $R_{m+1}$, then any join of $L_j$ with another loop will
encircle~$S_m$. In the case~$L_j$ does not encircle $R_{m+1}$,
loop~$L_j$ must be to the left of loop $L_{j'}$. Without loss of
generality, we may assume that the level of loop $L_j$ is not
bigger than the level of~$L_{j'}$. Then there are loops on the
level of loop $L_j$ to the right of $L_{j}$. In particular, $L_j$
gets joined with another loop~$L'$ which is to the right of~$L_j$.
Loop $L'$ must encircle at least one loop in the level $R_1,\dots,
R_t$. Moreover, since we assumed that~$L_j$ is the rightmost top
loop that encircles~$R_m$, loop~$L'$ must encircle at least
one~$R_{m'}$ with $m'\geq m+1$. So the join of $L_j$ and $L'$
encircles $S_m$.

To prove (2) repeat the argument for (1) with all indices $m$ and
${m+1}$ replaced by $\ell$ and $r$ and loop $S_m$ replaced by the
loop $S'$ whose left end coincides with the left end of $S_\ell$
and right end coincides with the right end of $S_{r-1}$. The whole
argument goes through and since loop $S'$ encircles all the loops
$S_\ell\dots, S_{r-1}$, the first part of the claim is proved.

For the second part of the claim notice that if $|J|+1=s$ then at
each step of the algorithm the number of top loops has to decrease
by exactly one. So the situation described in (2) can never
happen; if it did, the number of top loops would decrease by at
least $r-\ell$, since $S_\ell,\dots S_{r-1}$ are not top loops.

For $1\leq a\leq s-1$, let $L_{m_a}$ be the rightmost loop among
$L_1,\dots, L_s$ passing through column~$j_a$; in other words,
$m_a$ is the number of top loops of $\pp$ that pass through one of
the columns $j_1,\dots, j_a$. We claim that $m_a> a$.

It is clear that if we first collapse along one of the columns
$j_1,\dots, j_a$, then the number of top loops passing through one
of these columns decreases by one (since the total number of top
loops decreases by one at each collapse and all loops which do not
pass through these columns remain unchanged).

Let us show that the number of top loops passing through columns
$j_1,\dots, j_a$ does not increase if our first collapse is along
column $j_b$ with $b>a$. The only way this number might increase
is if one of the loops $L_m$ with $m> m_a$ is joined with another
loop $L$ which passes through $j_a$. If this happens, then
$L\subseteq L_{m-1}$. In particular, the top loop~$L_{m-1}$ passes
through column~$j_a$ and the top loop~$L_m$ does not. Hence
$m-1=m_a$. Moreover, the level of loop~$L_{m_a}$ is not smaller
than the level of loop~$L_{m_a+1}$, so that at least two of the
loops $L_1,\dots,L_{m_a+1}$ are in the highest level. Thus, after
collapsing, these loops produce at most~$m_a$ top loops passing
through one of the columns $j_1,\dots, j_a$.

Since at the end $\hat \pp_\II$ contains exactly one top loop
passing through columns $j_1,\dots, j_a$, we conclude, by
induction, that $m_a> a$. By the symmetric argument we conclude
that if $L_{n_a}$ is the leftmost loop passing through
column~$j_a$, then $n_a\leq a$. This proves that each $L_a$ passes
through columns~$j_{a-1}$ and~$j_a$. \hfill$\square$

\section{Lattices and their compatibility to Kostant pictures}
\label{sec:bases}

In this section $G=\gln$, so that $\Gr=\lgrgl$; we use the lattice
model and the terminology of Sections~\ref{sec:lattice-model},
\ref{sec:kostant-parameter} and~\ref{sec:MV-parametrization}. We
explain the basic relationship between Kostant pictures and
lattices, whereby loops correspond to basis vectors. Proofs of the
main results will depend on understanding the collection of
lattices ``compatible'' to a Kostant picture, with respect to
three possible degrees of compatibility.

\subsection{Bases of lattices} We say that a set of vectors
$y_1,\dots,y_k\in X$ is a \emph{basis} of a lattice~Y, if
$k=\dim_0(Y)$ and
\begin{eqnarray*}Y&=&Y_0\oplus
\langle y_1,\dots, y_k\rangle,
\end{eqnarray*}
where $\langle y_1,\dots, y_k\rangle$ is the $\CC$-span of
$y_1,\dots,y_k$.  Given a basis of $Y$, every vector~\mbox{$y\in
Y$} can be uniquely written as a linear combination of
$y_1,\dots,y_k$ plus a vector~$y_0\in Y_0$.

\subsection{$\pp$-flags}

Suppose $\pp$ is a Kostant picture and $Y$ is a lattice with basis
$\{y_L\}_{L\in \pp}$ indexed by the loops of $\pp$. Moreover,
assume that for each $L\in \pp$, we have $y_L \in V_L$ and $t
\cdot y_L \in Y_0 \oplus {\langle y_{L'} \rangle}_{L'\subset L}$.
The collection $\{Y_L \}_{L \in \pp}$ of $t$-invariant subspaces
of $Y$ defined by
\begin{eqnarray}
\label{eq:Y_L} Y_L&=& \big(Y_0\cap V_L\big) \oplus\langle
y_{L'}\rangle_{L'\subseteq L}
\end{eqnarray}
is called the $\pp$-\emph{flag} generated by the basis. Note that
if $L\subset L'$ then $Y_L \subset Y_{L'}$ and that $\dimz(Y_L)$
equals the number of $L' \in \pp$ with $L'\subseteq L$.

For example, the four vectors in the second picture of
Figure~\ref{fig:lattice-examples} form a basis that generates a
$\pp$-flag for the lattice, where $\pp$ is the second Kostant
picture of Figure~\ref{fig:Kostant-pictures}; similarly for the
third pictures.

By a \emph{subpicture} $\pp'$ of $\pp$ we mean a subset of the
loops of $\pp$; we say it is an {\em inner} subpicture if for
every loop $L$ of~$\pp'$, all the loops of $\pp$ encircled by $L$
are in $\pp'$. For a subpicture $\pp'$ let $V_{\pp'}$ be the sum
of all $V_L$ for loops $L$ in $\pp'$. A $\pp$-flag in $Y$
naturally extends to give a lattice $Y_{\pp'}$ in $V_{\pp'}$ for
every inner subpicture $\pp'$ of $\pp$: let $Y_{\pp'}$ be the sum
of the subspaces~$Y_L$ for $L$ in $\pp'$.
It is easy to see that $\dim_0(Y_{\pp'})$ equals the number of
loops in $\pp'$. In particular, we will use this construction for
the inner subpicture $\pp^L$ consisting of all loops encircled by
a fixed loop $L$ in $\pp$ and we let
$Y_L^\circ=Y_{\pp^L}+\big(Y_0\cap V_L\big)$. Note that $Y_L^\circ
\subseteq Y_L$ and $\dim(Y_L/Y_L^\circ)=1$.

\subsection{Lattices weakly compatible to Kostant pictures}
We say that a lattice $Y$ is \emph{weakly compatible} to $\pp$ if
there is a $\pp$-flag~$\{Y_L\}_{L\in\pp}$ in $Y$ such that for
every loop $L$ of $\pp$
\begin{romanlist}
\item
\label{property} {$\proj_{V_i}(Y_L)\neq \proj_{V_i}(Y_L^\circ)$
for the left and right ends $i$ of $L$,}
\end{romanlist}
where $\proj_V(W)$ stands for the orthogonal projection of $W$
onto $V$, treating the vectors $t^je_i$ as an orthogonal basis of
$X$.  Notice that by $t$-invariance, $\proj_{V_i}(Y_L^\circ) =
\langle t^{j}e_i, t^{j+1}e_i, t^{j+2}e_i, \dots \rangle$ and
$\proj_{V_i}(Y_L) = \langle t^{k}e_i, t^{k+1}e_i, t^{k+2}e_i,
\dots \rangle$ for some integers $j$ and $k$; condition
(\ref{property}) just says that $k=j-1$ rather than $k=j$.

For example, the second lattice in
Figure~\ref{fig:lattice-examples} is weakly compatible to the
second Kostant picture in Figure~\ref{fig:Kostant-pictures}. This
would not be the case, however, if, for instance, the vector
$t^{-2}e_1+3e_2+4e_3$ were changed to $3e_2+4e_3$. Similarly, the
third lattice is weakly compatible to the third Kostant picture.

\begin{prop}
\label{prop:weakly-compatible} A lattice $Y$ is weakly compatible
to a Kostant picture $\pp$ if and only if $\pp=\kst(Y)$, the
Kostant picture associated to $Y$.  Moreover, if $\pp=\kst(Y)$
then the $\pp$-flag $\{Y_L\}_{L\in\pp}$ in $Y$ is unique.
\end{prop}

\begin{proof} Assume we are given a $\pp$-flag
$\{Y_L\}_{L\in\pp}$ satisfying (\ref{property}). Let
$L_1\subset\dots\subset L_k$ be all the loops of~$\pp$ with left
end at column~$\ell$ and right end at column $r$.
Then, by definition, the number of such loops~$L$ in~$\kst(Y)$
is equal to
\begin{align*}
n_L
&=\dim\Big(\big(Y\cap V_L\big)\Big/\big((Y\cap V_{L_+})+(Y\cap
V_{L_-})\big)\Big)=\dim (Y_{L_k}/Y^\circ_{L_1})=k,
\end{align*}
where we used condition (\ref{property}) in the last two
equalities. Hence~$\kst(Y)=\pp$.

Conversely, assume that $Y$ is a lattice with $\kst(Y)=\pp$. We
inductively define a $\pp$-flag in $Y$.  For a loop $L$ of $\pp$
that encircles no loops, let
\begin{equation*}
Y_L=Y\cap t^{-1}\big(Y_0\cap V_L\big).
\end{equation*}
By induction assume the spaces $Y_{L'}$ are defined for all loops
$L'$ encircled by loop $L$ and define a $\pp^L$-flag in
$Y_L^\circ$, where the loops in $\pp^L$ satisfy condition
(\ref{property}). So, by hypothesis, we must have
\begin{equation}
\label{eq:1} Y\cap V_{L_+}\subset Y^\circ_L\ \ \text{and}\ \ Y\cap
V_{L_-}\subset Y^\circ_L.
\end{equation}
Set
\begin{equation}
\label{eq:2} Y_L=Y\cap t^{-1}Y^\circ_L.
\end{equation}

Because of (\ref{eq:1}) no vector in $Y_L-Y^\circ_L$ can restrict
trivially to $V_i$ for the left or right end $i$ of $L$; since
$\dim\Big(\proj_{V_i}(Y_L)/ \proj_{V_i}(Y_L^\circ)\Big)\leq 1$ for
the left and right ends, we conclude that $\dim
(Y_L/Y^\circ_L)\leq 1$.

To show $\{Y_L\}_{L\in\pp}$ is a $\pp$-flag it remains to show
$Y_L\neq Y^\circ_L$. Assume $L$ passes through columns
$\ell,\dots,r$. Let $L_1,\dots,L_k$ be all the loops of $\pp$ with
left end $\ell$ and right end $r$. So $L=L_m$ for some $m$. If
$Y_L= Y^\circ_L$, then $Y^\circ_L = Y\cap t^{-1}Y^\circ_L=Y\cap
V_L$. Hence $n_L=\dim\Big((Y\cap V_L)/\big((Y\cap V_{L_+})+(Y\cap
V_{L_-})\big)\Big)=m-1<k$, a contradiction.

This finishes the construction of the collection
$\{Y_L\}_{L\in\pp}$.  Its uniqueness follows immediately from the
fact that for any $\pp$-flag we have $Y_L\subset Y\cap
t^{-1}Y^\circ_L$ and hence definition~(\ref{eq:2}) is the only one
possible.
\end{proof}

We have just shown that all lattices in $\wcomp(\pp,\lambda)$ are
weakly  compatible to $\pp$.
We now identify $\wcomp(\pp,\lambda)$ with an open Zariski subset
of a product of projective spaces.
For a loop $L$ with left end $\ell$ and right end $r$,
let $\PP^L$ be the projective space of complex lines inside
$\langle e_\ell,\dots,e_r\rangle$. Set
$\PP^\pp=\prod_{L\in\pp}\PP^L$.
If loop $L$ encircles no loop with right end $r$ then let $m=r$;
otherwise let $M$ be the largest loop with right end $r$ encircled
by $L$, and let $m$ be its left end; note that $m=\ell$ is
possible. Let $\tilde \PP^L$ be the set of lines in $\langle
e_\ell,\dots,e_r\rangle$ that project nontrivially onto $\langle
e_\ell\rangle$ and~$\langle e_{m}\rangle$.
Define~$\tilde\PP^\pp=\prod_{L\in\pp}\tilde\PP^L$.

\begin{prop}
\label{prop:M-is-smooth} $\wcomp(\pp,\lambda)$ is isomorphic to
$\tilde \PP^\pp$.
\end{prop}

\begin{proof}
Fix $\pp$ and $\lambda$. Without loss of generality, we may
assume, using a shift by an appropriate coweight, that $\lambda$
is such that $Y_0=X_0$ for any lattice $Y\in\wcomp(\pp,\lambda)$.
We will construct an algebraic map $\pi: \tilde\PP^\pp \to \lgr$
which will provide the required isomorphism with the image
$\wcomp(\pp,\lambda)$.
To do so, given a point $p\in \tilde\PP^\pp$, we will construct a
basis $\{y_L\}_{L\in\pp}$ that generates a
\mbox{$\pp$-flag}~$\{Y_L\}_{L\in\pp}$ satisfying~(\ref{property}).
Moreover, each $y_{L}$ will be of the form $t^je_\ell+z$ for some
$j$ and $z\in V_{\ell+1}\oplus \dots \oplus V_{r}$, where $\ell$
and $r$ are the left and right ends of~$L$.

The construction proceeds by induction on the number of loops of
$\pp$. If $\pp$ contains only one loop~$L$, say with left end
$\ell$ and right end $r$, then $p$ defines a line in $\langle
e_\ell,\dots,e_r\rangle$ with nontrivial projections onto $\langle
e_\ell \rangle$ and $\langle e_r\rangle$. Let $x$ be the unique
point in the line $p$ with $\proj_{\langle
e_\ell\rangle}x=e_\ell$.
Let $Y_L=Y_0\oplus t^{-1}p$ and $y_L=t^{-1}x$. Clearly, the
constructed lattice $Y=Y_L$ is in~$\wcomp(\pp,\lambda)$.

Now assume that a loop $L$ is not encircled by any other loop of
$\pp$. Then $\pp'=\pp-\{L\}$ is an inner subpicture of $\pp$. By
inductive assumption, given~$p\in \tilde \PP^\pp$ we can construct
a $\pp'$-flag $\{Y_{L'}\}_{L'\in\pp'}$ for $Y_{\pp'}$ together
with a basis $\{y_{L'}\}_{L'\in\pp'}$ for which
each $y_{L'}$ is of the form $t^j e_{\ell'}+z$ for some $j$ and
$z\in V_{\ell'+1}\oplus \dots \oplus V_{r}$ with $\proj_{V_{r'}}
y_{L'}$ nonzero; here $\ell'$ and $r'$ are the left and right ends
of $L'$. We need to show how to construct $Y_L$ and $y_L$, given
the line $p_L$ in $\langle e_\ell,\dots e_r\rangle$, the
$L^{\rm{th}}$ component of $p$.

Recall that $Y_L$ must lie inside $t^{-1}Y^\circ_L$ and $\dim
(Y_L/Y^\circ_L)=1$.  So, if we can identify
$t^{-1}Y^\circ_L/Y^\circ_L$ with $\langle e_\ell,\dots e_r\rangle$
then the line $p_L$ defines a line inside
$t^{-1}Y^\circ_L/Y^\circ_L$ and hence defines $Y_L$.

For $\ell\leq k\leq r$, let $L_k$ be the largest loop encircled by
$L$ whose left end is at column $k$ and set $y_k=y_{L_k}$; if no
loop encircled by $L$ has column $k$ as its left end, set
$y_k=e_k$. We claim that $y_\ell,\dots y_r$ is a $\CC[t]$-basis
of~$Y^\circ_L$; that is, $Y^\circ_L$ is the direct sum of $\langle
y_k\rangle_t$, the $\CC[t]$-span of the vector $y_k$.

To prove it, let us show that any vector $y\in Y^\circ_L$ can be
uniquely written as a $\CC[t]$-linear combination of the vectors
$y_k$. Let $k$ be the smallest number such that $\proj_{V_k} y\neq
0$. Then, since $y_k$ has the form $t^j e_k + z$, there is a
unique linear combination $\sum_j c_jt^jy_k$ such that its
difference with $y$ projects trivially onto $V_k$. If we can show
that~$c_j=0$ for $j<0$, then induction on $k$ proves that
$y_\ell,\dots y_r$ is a $\CC[t]$-basis of $Y^\circ_L$. Since $y$
can be written uniquely as a linear combination of basis vectors
plus a vector from $Y_0$, it is enough to assume $y$ is a basis
vector $y_{L'}$ for a loop $L'$ whose left end is at column~$k$. But
since $L'\subseteq L_k$, every~$c_j$ must be~$0$ for negative~$j$.

This implies that that we can identify $t^{-1}Y^\circ_L/Y^\circ_L$
with $\langle t^{-1}y_\ell,\dots t^{-1}y_r\rangle$. Mapping
$t^{-1}y_k$ to $e_k$ provides the isomorphism between
$t^{-1}Y^\circ_L/Y^\circ_L$ and $\langle e_\ell,\dots,e_r\rangle$.
So, as explained before,  the line $p_L$ in $\langle
e_\ell,\dots,e_r\rangle$ defines $Y_L$. To define the vector
$y_L$, let $x_L=e_\ell + \sum_{k=\ell+1}^r a_k e_k$ be the unique
point in $p_L$ with $\proj_{\langle e_\ell\rangle}x_L =e_\ell$.
Then set $y_L=t^{-1}y_\ell+ \sum_{k=\ell+1}^r a_k t^{-1}y_k$.

To finish the construction of $\pi$, we have to show
that~(\ref{property}) holds for loop $L$. Since
$\proj_{V_\ell}t^{-1}y_\ell\notin \proj_{V_\ell}Y^\circ_L$ we
conclude $\proj_{V_\ell} y_L\notin \proj_{V_\ell}Y^\circ_L$; hence
(\ref{property}) holds for the left end of $L$. If $L$ encircles
no loop with right end $r$ then $\proj_{V_r} y_k = 0 $ for $k<r$
so that $\proj_{V_r} y_L =a_r t^{-1} y_r = a_r t^{-1} e_r \notin
\proj_{V_r} Y^\circ_L$. Otherwise, recall that $M$ is the largest
loop encircled by $L$ with right end $r$; it has left end $m$,
and~$p_L$ projects nontrivially onto $\langle e_m \rangle$. Then
$M=L_m$ and $a_m\neq 0$.
Since $\proj_{V_r}t^{-1}y_{m}\notin \proj_{V_r}Y^\circ_L$  we
conclude $\proj_{V_r} y_L\notin \proj_{V_r}Y^\circ_L$; hence
(\ref{property}) holds for the right end of $L$ in either case.

It remains to show that for every $Y\in \wcomp(\pp,\lambda)$ there
is a unique $p\in \tilde \PP^\pp$ with $\pi(p)=Y$. Let
$\{Y_L\}_{L\in\pp}$ be the $\pp$-flag for $Y$. Then from the above
discussion it is clear that we can always construct a basis
$\{y_L\}_{L\in\pp}$ such that for each $L$,
and basis vectors $y_\ell, \dots, y_r$ as before, we have
$y_L=t^{-1}y_\ell+ \sum_{k=\ell+1}^r a_k t^{-1}y_k$. Given such a
basis of $Y$, define $p\in\tilde \PP^L$ by setting $p_L$ to be the
line passing through $e_\ell+\sum_{k=\ell+1}^r a_k e_k$. Clearly,
$\pi(p)=Y$ and such $p$ is unique.
\end{proof}

\subsection{Lattices compatible to Kostant pictures}
We say that a lattice $Y$ is \emph{compatible} to $\pp$ if it has
a $\pp$-flag $\{Y_L\}_{L\in\pp}$ such that for every loop $L$ of
$\pp$,
\begin{romanlist'}
\item
\label{property'} {$\proj_{V_i}(Y_L)\neq \proj_{V_i}(Y_L^\circ)$
for every column $i$ that $L$ passes through.}
\end{romanlist'}
Since (\ref{property'}$'$) is stronger than (\ref{property}),
every lattice compatible to $\pp$ is weakly compatible to $\pp$.
We denote by $\comp(\pp,\lambda)$ the set of all lattices in
$\wcomp(\pp,\lambda)$ that are compatible to $\pp$.

For example, the third lattice in
Figure~\ref{fig:lattice-examples} is compatible to the third
Kostant picture in Figure~\ref{fig:Kostant-pictures}.  But the
second lattice in Figure~\ref{fig:lattice-examples} is not
compatible to the second Kostant picture in
Figure~\ref{fig:Kostant-pictures}, since (\ref{property'}$'$)
fails for column~4 and the loop encircling dots~$3$ and~$4$; they
are, however, weakly compatible.

\begin{prop}
\label{prop:compatible}
$\comp(\pp,\lambda)$ is a dense Zariski open subset of
$\wcomp(\pp,\lambda)$.
\end{prop}

\begin{proof}
For an inner subpicture $\pp'$ of $\pp$, denote by
$\comp_{\pp'}(\pp,\lambda)$ the set of all lattices $Y\in
\wcomp(\pp,\lambda)$ whose $\pp$-flag
satisfies~(\ref{property'}$'$) for all loops in $\pp'$. We will
show that if $L\in\pp$ is a loop for which $\pp'\cup\{L\}$ is also
an inner subpicture then $\comp_{\pp'\cup L}(\pp,\lambda)$ is a
dense Zariski open dense of $\comp_{\pp'}(\pp,\lambda)$.  The
proof then follows by induction on the number of loops of $\pp$.

Given $Y\in \comp_{\pp'}(\pp,\lambda)$, consider the basis and
$\pp$-flag constructed for it in the proof of
proposition~\ref{prop:M-is-smooth}. For loop $L$, there was a
$\CC[t]$-basis $y_\ell,\dots,y_r$ of $Y^\circ_L$, and $y_L$ was
defined as a linear combination $\sum a_i t^{-1}y_i$ with
$a_\ell=1$. For $\ell\leq k\leq r$, let $j_k$ be defined by
$\proj_{V_k} Y^\circ_L=\langle t^{j_k}e_k\rangle_t=\langle
t^{j_k}e_k,t^{{j_k}+1}e_k,\dots\rangle$. Then
condition~(\ref{property'}$'$) for loop $L$ and column $k$ is
equivalent to $\sum_i a_i \proj_{\langle t^{j_k}e_k\rangle}
y_i\neq 0$. Combining such conditions for all columns $\ell$
through $r$ defines a dense Zariski open subset of
$\comp_{\pp'}(\pp,\lambda)$; this subset is $\comp_{\pp'\cup
L}(\pp,\lambda)$.
\end{proof}

\subsection{Lattices strongly compatible to Kostant pictures}
Recall that in Section~\ref{algorithm} we defined the collapse
$\hat \pp_i$ of a Kostant picture $\pp$ along column $i$. We also
define the collapse $\hat Y_i$ of a lattice $Y$ along column $i$
as the intersection of $Y$ with the direct sum of all columns
other than $i$.  Just as for~$\hat \pp_i$, we may view $\hat Y_i$
as a lattice one rank lower by removing the collapsed column and
renumbering the remaining columns.

\begin{lemma}
\label{lem:induction} If $Y$ is compatible to $\pp$ then $\hat
Y_i$ is weakly compatible to $\hat\pp_i$.
\end{lemma}

\begin{proof} Let  $\{Y_L\}_{L\in\pp}$
be the $\pp$-flag for $Y$ satisfying (\ref{property'}$'$). We must
produce a $\hat \pp_i$-flag~$\{\hat Y_M\}_{M\in\hat \pp_i}$ for
the lattice $\hat Y_i$ satisfying (i).

As explained in Section~\ref{algorithm}, every loop $M$ of $\hat
\pp_i$ is either a loop of $\pp$ or a join  of two loops of $\pp$.
If $M$ is a loop of~$\pp$, set $\hat Y_M=Y_M$. If $M$ is the join
of loops $L'$ and~$L''$, set $\hat Y_{M}= \hat V_i \cap (Y_{L'} +
Y_{L''})$, where $\hat V_i=\bigoplus_{j\neq i} V_j$.

Choose a basis $\{y_L\}_{L\in\pp}$ of $Y$ that generates the
$\pp$-flag $\{Y_L\}_{L\in\pp}$.
For a loop $L$ passing through column $i$, denote by $j_L$ the
largest power of~$t$ for which~$y_L$ projects nontrivially onto
$\langle t^{-j}e_i\rangle$. Without loss of generality we can
assume that $\proj_{V_i} y_L= t^{-j_L}e_i$. Indeed, this can be
achieved by adding to each $y_L$ a linear combination of the
vectors $t^jy_L$ for $j\geq 1$, multiplying~$y_L$ by a scalar, and
adding a vector from $Y_0$.

Recall that in the definition of $\hat \pp_i$,
the loops of $\pp$ passing through column~$i$ were separated into
levels, such that no loop encircles a loop from a higher level.
Since $\{Y_L\}_{L\in\pp}$ satisfies (\ref{property'}$'$), it is
easy to see that $j_L$ is equal to $\base_i(Y)$ plus the level of
loop $L$. In particular, if loop $M$ of~$\hat \pp_i$ is the join
of two loops $L'$ and $L''$, necessarily from the same level,
then $j_{L'}=j_{L''}$. Hence we can define $\hat
y_{M}=y_{L'}-y_{L''}\in \hat Y_i$. For a loop $M$ of $\hat \pp_i$
that is a loop of $\pp$, set $\hat y_{M}=y_M$.

Clearly the basis $\{\hat y_{M}\}_{M\in\hat \pp_i}$ generates the
$\hat \pp_i$-flag $\{\hat Y_{M}\}_{M\in\hat \pp_i}$ for $\hat
Y_i$. It remains to show that for every loop $M$ of $\hat \pp_i$,
property~(\ref{property}) holds. This is obvious if $M$ is a loop
of $\pp$. If $M$ is the join of~$L'$ and~$L''$, then the left end
$\ell$ of $M$ is the left end of $L'$, and $L''$ does not pass
through this column. Hence
$$
\proj_{V_\ell}{\hat
y_{M}}=\proj_{V_\ell}(y_{L'}-y_{L''})=\proj_{V_\ell}(y_{L'})
\notin \proj_{V_\ell} (Y^\circ_{L'}) = \proj_{V_\ell} (\hat
Y^\circ_{M}).
$$
This proves~(\ref{property}) for the left end, and an analogous
argument proves it for the right end.
\end{proof}

The definition of \emph{strongly compatible} is inductive on the
number of columns $n$: For $n=2$, a~lattice~$Y$ is said to be
strongly compatible to a Kostant picture $\pp$ if it is compatible
to $\pp$. For $n>2$ a lattice $Y$ is said to be strongly
compatible to $\pp$ if it is compatible to $\pp$ and if for every
column $i$ the lattice $\hat Y_i$ is strongly compatible to $\hat
\pp_i$. We denote by $\scomp(\pp,\lambda)$ the set of all lattices
in $\comp(\pp,\lambda)$ that are strongly compatible to $\pp$.

For example, the third lattice in
Figure~\ref{fig:lattice-examples} is compatible but not strongly
compatible to the third Kostant picture in
Figure~\ref{fig:Kostant-pictures}. Indeed, if we collapse the
Kostant picture along column~$4$, the resulting Kostant picture
contains only two loops: the rightmost loop of length~$1$, which
is unchanged, and a new loop~$M$.  The lattice $\hat Y_3$ is
generated by two corresponding basis vectors: $te_5+t^{-1}e_6$ and
$t^2e_1+t^2e_2-te_5$ (where we have not renumbered the columns
after the collapse). Since the latter vector contains no $e_3$
component, we see that (\ref{property'}$'$) fails for loop~$M$ and
column $3$. If, however, we alter the original vector
$te_2+2te_3+3te_4+te_5$ by adding to it any nonzero multiple
of~$t^{2}e_3$, then the lattice becomes strongly compatible to the
Kostant picture.

\begin{remark}
We will see later, in Theorem~\ref{thm:strong-moment}, that
$\scomp(\pp,\lambda)$ is the largest piece of the decomposition
into torus orbit types contained in $\wcomp(\pp,\lambda)$; that
is, the strongly compatible lattices are precisely those for which
the moment map image of the closure of the torus orbit through the
lattice is the entire MV-polytope.
\end{remark}

We need notation for collapsing lattices along two or more
columns. For a set $I=\{i_1,\dots,i_k\}\subseteq\{1,\dots,n\}$
define
$$
V_I=V_{i_1}\oplus \dots \oplus V_{i_k}.
$$
In particular, if $I=[\ell..r]=\{\ell,\dots,r\}$ and $L$ is the
loop with left end $\ell$ and right end $r$, then $V_L=V_I$. For a
lattice $Y$, set $Y_I=Y\cap V_I$. Denote by $\hat I$ the
complement of $I$ in $[1..n]$. Let $\hat Y_I = Y_{\hat I}$, the
collapse of the lattice $Y$ along columns $I$; in particular when
$I=\{i\}$ we get $\hat Y_i$.

\def\II{\mathbb I}

For collapsing a Kostant picture along several columns $I$, the
order in which the columns collapse matters, a priori.  Let
$\II=(i_1,\dots,i_k)$ be an ordering of $I$.
As in Section~\ref{algorithm} denote by $\hat \pp_{\II}$ the
Kostant picture produced out of $\pp$ by collapsing columns $\II$
in the given order.

\begin{remark}
Once we have proved, in Proposition~\ref{strongly-compatible},
that there always exists a lattice $Y$ strongly compatible to
$\pp$, we will know that the order of collapse doesn't matter;
indeed, $\hat \pp_{\II}=\pp(\hat Y_I)$ for any ordering of the set
$I$. Thus, we will have proved Theorem~\ref{thm:commutativity}
from Section~\ref{sec:combinatorics}.
\end{remark}

Notice that a lattice $Y$ is strongly compatible to $\pp$ if and
only if $\hat Y_I$ is compatible to $\hat\pp_{\II}$ for every
ordered set $\II$. Indeed, if $Y$ is  strongly compatible to
$\pp$, then every $\hat Y_I$ must be strongly compatible (hence
compatible) to $\hat \pp_{\II}$. Conversely, if the lattice $\hat
Y_I$ is compatible to $\hat \pp_{\II}$ for every $\II$, then, by
induction on $n$, we can assume every $\hat Y_i$ is strongly
compatible to $\hat \pp_i$; therefore $Y$ is strongly compatible
to~$\pp$.

\begin{prop}
\label{strongly-compatible} $\scomp(\pp,\lambda)$ is a dense
Zariski open subset of $\comp(\pp,\lambda)$.
\end{prop}

\begin{proof}
We fix $\pp$ and $\lambda$ throughout the proof.  For an ordered
set $\II=(i_1,\dots,i_k)$, let $\II_m=(i_1,\dots,i_m)$ for $1\leq
m\leq k$.  Let $I$ and $I_m$ denote the corresponding unordered
sets.
Let~$S_\II$ be the subset of $\comp(\pp,\lambda)$ of lattices $Y$
for which $\hat Y_{I_m}$ is compatible to $\hat\pp_{\II_m}$ for
every $1\leq m\leq k$. Then, as mentioned before,
$\scomp(\pp,\lambda)$ is the intersection of all $S_\II$. Since
every~$S_\II$ is clearly an open Zariski subset of
$\comp(\pp,\lambda)$, it is enough to show that every $S_\II$ is
not empty.

Fix $\II$. Let $M^m$ be the loop of $\hat\pp_{\II_m}$ whose left
end is furthest left; if there is more than one such loop, let it
be the largest. Similarly, let $M$ be the largest leftmost loop of
$\pp$. Then, recalling terminology from
Section~\ref{sec:combinatorics}, if a loop of $\hat\pp_{\II_m}$ is
a join of loops, one of which is~$M$, this loop must be~$M^m$, and
we say that the loop $M$ survives until step $m$; that is, $M^m$
is the only loop of $\hat\pp_{\II_m}$ whose ancestry contains $M$.

Set $\qq=\pp-\{M\}$ and let $Q=Y_\qq+Y_0$ denote the lattice
associated to this inner subpicture. It is easy to see that if $M$
does not survive until step $m$ then $\hat \qq_{\II_m} =
\hat\pp_{\II_m}$, and that if $M$ survives until step $m$ then
$\hat \qq_{\II_m} = \hat\pp_{\II_m}-\{M^m\}$.

We will use induction on the number of loops of $\pp$, and, within
that, induction on the length $k$ of $\II$ to prove:
\renewcommand{\labelitemi}{($\star$)}
\begin{itemize} \item
There exists a lattice $Y\in S_\II$ for which $Q$ is strongly
compatible to $\qq$.
\end{itemize}
The base case for the outer induction is for a Kostant picture for
no loops, which is trivial.  The base case for the inner induction
is that there exists a lattice $Y$ for which $Q$ is strongly
compatible to $\qq$.  Since $\qq$ contains one fewer loop, we
inductively know the existence of a lattice $Q$ strongly
compatible to $\qq$, which may be extended by a single vector
$y_M$ to give $Y$.

For the induction on $k$, pick a lattice $Y$ that satisfies
($\star$) for $\II_{k-1}$. First consider the case that~$M$ does
not survive until step $k$. We claim that the same lattice $Y$
satisfies ($\star$) for~$\II_k$. It is enough to show that $\hat
Y_{I_k}=\hat Q_{I_k}$ since by inductive assumption $\hat Q_{I_k}$
is compatible to $\hat \pp_{\II}= \hat \qq_{\II}$. Obviously $\hat
Q_{I_k} \subseteq \hat Y_{I_k}$. Since $\hat Y_{I_{k-1}}$ is
compatible to $\hat \pp_{\II_{k-1}}$ we know that $\hat Y_{I_{k}}$
is weakly compatible to $\hat \pp_{\II_{k}}$ by
Lemma~\ref{lem:induction}. Hence $\dimz(\hat Y_{I_k})=|\hat
\pp_{\II_k}|= |\hat \qq_{\II_k}|=\dimz(\hat Q_{I_k})$, and since
$Y_0=Q_0$, we must have $\hat Y_{I_k}=\hat Q_{I_k}$.

If, however, $M$ survives until step $k$, we will show that there
exists a lattice arbitrarily close to~$Y$ that satisfies ($\star$)
for $\II_k$. Let $\{y_L\}_{L\in \pp}$ be a basis of $Y$ that
generates its $\pp$-flag.  To shorten some notations, let us
introduce the extended Kostant picture $\tilde \pp=\pp\cup
[1..n]$. So $L\in \tilde\pp$ means either $L$ is a loop of $\pp$
or an integer between $1$ and $n$. (It might help to think about
the integer $i$ as a loop of length zero passing through the
single column~$i$.) Of course we say that $L$ encircles
integer~$i$ if and only if $L$ passes through column~$i$. Set
$y_i=t^{-\base_i (Y)} e_i$.

Since for a loop $L$ of $\pp$ we have $t\cdot y_L\in Y^\circ_L$,
we can uniquely write the corresponding basis vector~as
\begin{equation}
\label{eq:linear} y_L\eqz \ t^{-1}\sum_{L'\in \tilde \pp}
a_{L}^{L'} y_{L'},
\end{equation}
where $v\eqz\ w$ means $v-w\in Y_0$, and $a_L^{L'}=0$ unless
$L'\subset L$. Conversely, any set of numbers $a_L^{L'}$ with this
property uniquely defines a lattice. We will perturb $Y$ by
specifying how to change these numbers. As long as the changes we
make are arbitrarily small, we may assume that the conditions in
($\star$) still hold since they are open conditions. So to show
that ($\star$) holds for $\II$, we need only find an arbitrarily
small perturbation of numbers $a_L^{L'}$ such that $\hat Y_{I}$ is
compatible to $\hat \pp_\II$.

Let $\{ Y_L \}_{L\in \hat \pp_\II}$ be a $\hat \pp_{\II}$-flag in
$\hat Y_I$. The columns of $\hat Y_I$ are naturally identified
with those $i$ not in~$\II$, and we will keep this numbering. We
need only show that $Y_{M^k}$ satisfies (\ref{property'}$'$),
since the other $Y_L$ satisfy it by the second part of the
induction hypothesis. By Lemma~\ref{lem:induction}, $\hat Y_{I}$
is weakly compatible to~$\hat \pp_\II$.  So there exists a vector
$y\in \hat Y_I$ such that $Y_{M^k}=Y^\circ_{M^k}\oplus \langle y
\rangle$.
If $y$ satisfies $\proj_{V_i} y\notin \proj_{V_i} Y_{M^k}^\circ$
for every $i\notin \II$ that ${M^k}$ passes through, then we are
done. If not there exists a vector $v\in Y_{M^k}^\circ$, as small
as we would like it to be, such that $y+t^{-1}v$ satisfies this.
Suppose $M^k$ is the join of loops $M=L_1,\dots, L_s$. We will
show how to change the numbers $a_{L_j}^{L'}$ by small amounts so
that $y$ is replaced by $y+t^{-1}v$.

Let $\mm$ denote the Kostant subpicture of $\tilde \pp$ consisting
of $L_1,\dots,L_s$ together with all the loops they encircle, and
let $\mm^\circ = \mm - \{ L_1, \dots, L_s \}$. So $Y_\mm =
Y_{L_1}+\dots +Y_{L_s}$ and $Y_{\mm^\circ}= Y^\circ_{L_1}+\dots
+Y^\circ_{L_s}$. By the proof of Lemma~\ref{lem:induction} we know
that $Y_{M^k}=V_I \cap Y_\mm$ and $Y^\circ_{M^k}=V_I\cap
Y_{\mm^\circ}$. So we can uniquely write
\begin{align}
\label{eq:linear1} y&\eqz\sum_{L\in \mm } b_L y_L, \\
\label{eq:linear2} v&\eqz\sum_{L\in \mm^\circ } c_L y_L.
\end{align}
If $b_{L_j}\neq 0$ for all $j=1,\dots,s$, the perturbation is
easy: for every nonzero $c_L$, pick a loop $L_j$ that encircles
$L$ and change $a_{L_j}^L$ to $a_{L_j}^L +\frac{c_L}{b_{L_j}}$.
This changes $y_{L_1},\dots,y_{L_s}$ so that the same linear
combination~(\ref{eq:linear2}) now equals $y+t^{-1}v$; so $\hat
Y_I$ is now compatible to $\hat \pp_\II$.

It remains to show that we can perturb $Y$ slightly so that every
$b_{L_j}$ is nonzero. Let $J=(j_1,\dots,j_{|J|})\subseteq I$
denote those columns of $I$ passed through by at least one of the
loops $L_1,\dots L_s$.
Let us show that $J$ contains exactly $s-1$ elements, which will
prove Proposition~\ref{prop:s-1}.

Let $\ell$ be the left end of $L_1$ and $r$ the right end of
$L_s$, and let $\bar J$ be the complement of $J$ in $[\ell..r]$.
Consider these finite-dimensional vector spaces:
\begin{align*}
&U=t^{-1}Y_{\mm^\circ}/Y_{\mm^\circ}, \\
&V = Y_{\mm} / Y_{\mm^\circ}, \\
&Z=\big(Y_{M^k}+Y_{\mm^\circ}\big)/ Y_{\mm^\circ}, \\
&W= \big((t^{-1}Y_{\mm^\circ} \cap V_{\bar J}) +
Y_{\mm^\circ}\big)/ Y_{\mm^\circ}.
\end{align*}
Then $V$, $Z$ and $W$ are subspaces of $U$. Moreover, $\dim V=s$,
$\dim Z=1$, $\dim W=|\bar J|$, and $\dim U=r-\ell+1=|J|+|\bar J|$,
since for any lattice $Y^\circ$, the dimension of
$t^{-1}Y^\circ/Y^\circ$ equals the number of columns.
Since $V\cap W=Z$, we have $\dim V +\dim W \leq \dim U+1$, so that
$s-1=\dim V-1\leq \dim U -\dim W=|J|$. Together with
Claim~\ref{claim2} we conclude $|J|=s-1$. So
$J=(j_1,j_2,\dots,j_{s-1})$ and, by the second part of the claim,
each $L_m$ passes through columns $j_{m-1}$ and $j_m$.

Now we will work in the $s-1$~dimensional vector space $U/W$. Let
$\tilde y_{L_j}\in U/W$ be the vectors corresponding to $y_{L_j}$.
Then, since $y\in W$, we clearly have
$$
0=\sum_{j=1}^s b_{L_j} \tilde y_{L_j}.
$$
Also let $\tilde y_j\in U/W$ correspond to $t^{-1} y_j$; note that
if $j\in J$ then $\tilde y_j$ is nonzero.

Consider perturbations just of the numbers $a^i_{L_i}$ and
$a^{i-1}_{L_{i}}$, which are allowed by the second part of
Claim~\ref{claim2}. We claim that there exists an arbitrarily small
change of these numbers such that every  $s-1$ vectors among
$\tilde y_{L_1},\dots,\tilde y_{L_{s}}$ are linearly independent
and, in particular, all $b_{L_j}$ are nonzero.

First, let's perturb these numbers to make $\tilde
y_{L_1},\dots,\tilde y_{L_{s-1}}$ independent. We will prove by
induction on $m$ that the vectors $\tilde y_{L_1},\dots,\tilde
y_{L_{m}}$ can be made independent in~$U/W_m$, where
$$
W_m= \big((t^{-1}Y_{\mm^\circ} \cap V_{\bar J_m}) +
Y_{\mm^\circ}\big)/ Y_{\mm^\circ},
$$
and $J_m=(j_1<\dots<j_m)$. Indeed, if this is true, then the
vectors $\tilde y_{L_1},\dots,\tilde y_{L_{m}}, \tilde
y_{j_{m+1}}$ are linearly independent inside $U/W_{m+1}$.
Moreover, if $\tilde y_{L_1},\dots,\tilde y_{L_{m+1}}$ are
linearly dependent inside $W_{m+1}$, the dependency relation must
have a nonzero coefficient in front of $\tilde y_{L_{m+1}}$;
therefore, an alteration of $a^{m+1}_{L_{m+1}}$ by any amount,
which can be taken arbitrarily small,
will make them independent.

For the $s-1$ vectors $\tilde y_{L_1},\dots,\tilde
y_{L_{i-1}}\tilde y_{L_{i+1}},\dots,\tilde y_{L_{s}}$, the same
argument works, except that we have to alter numbers
$a^1_{L_1},\dots,a^{i-1}_{L_{i-1}},a^{i}_{L_{i+1}},\dots,a^{s-1}_{L_{s}}$.
Obviously we can successively choose
$a^1_{L_1},\dots,a^{s-1}_{L_{s}}$
arbitrarily small so that linear independence
simultaneously holds for each of these $s$~subsets of
$s-1$~vectors.
\end{proof}

\section{Decompositions of $\Gr$ and the moment map.}
\label{sec:decompositions} Throughout this section, $G$ is any
connected simply-connected semisimple complex algebraic group, and
$\Gr$ its loop Grassmannian. This section contains the facts about
$\Gr$ that we can prove for all types. We review three known
decompositions of $\Gr$ and discuss the moment map images of the
pieces. The main theorem states that all moment map images of
compact irreducible algebraic subvarieties of $\Gr$ have the same
shape: an intersection of certain cones spanned by coroots.

\subsection{Decompositions of $\Gr$}

First, the left action of $G(\OO)$ on $\Gr$ decomposes it into
orbits.  We let~$G_\lambda$ denote the $G(\OO)$-orbit through
$\underline{\lambda}$. The coweight $\lambda$ is determined up to
the action of the Weyl group, and we have
\begin{eqnarray*}
\Gr&=&\bigcup_{\lambda\in \Lambda^+}G_\lambda
\end{eqnarray*}
where $\Lambda^+$ is the set of all dominant coweights. Each
closure $\overline{G_\lambda}$ is a finite-dimensional projective
variety and is the union of the $G_\mu$ with $\mu\leq\lambda$ and
$\mu\in \Lambda^+$. (Here $\leq$ is the usual partial order on
coweights: $\mu\leq\lambda$ means that $\lambda-\mu$ is a
nonnegative integral linear combination of simple coroots.)

Second, we have the finer \emph{Bruhat decomposition} of $\Gr$
into $B(\OO)$-orbits, where $B$ is a Borel subgroup of $G$. Each
orbit contains exactly one torus fixed point $\underline{\lambda}$
and is an affine cell~\cite{AP}.
\begin{eqnarray*}
\Gr&=&\bigcup_{\lambda\in \Lambda}B_\lambda
\end{eqnarray*}

Third, we need the decomposition of $\Gr$ into $\unip(\KK)$-orbits,
where $\unip$ is the unipotent radical of $B$.  Each orbit is
infinite-dimensional and again contains exactly one $\underline{\lambda}$;
we denote it by $S_\lambda$.
\begin{eqnarray*}
\Gr&=&\bigcup_{\lambda\in \Lambda}S_\lambda
\end{eqnarray*}
The closure $\overline{S_\lambda}$ is the union
$\bigcup_{\mu\leq\lambda} S_\mu$.  Note that we get such a
decomposition of $\Gr$ for each element~$w$ of the Weyl group: we
set $S_\lambda^w=wN(\KK)w^{-1}\underline\lambda$. The following
lemma uses the torus action to identify in which piece of this
decomposition a given $x \in \Gr$ lies.

Every coweight $\beta$ defines a one-parameter subgroup of the
maximal torus~$T$ whose elements are given by $\exp(s\beta)$
for~$s\in \CC$. If $\tau=e^s\in \CC-\{0\}$
we write $\tau^\beta= \exp(s\beta)$.

\begin{lemma}
\label{lem:limit} A point $x\in \Gr$ is in $S_\lambda^w$ if and
only if
\begin{equation*}
\lim_{\tau\to 0} w\tau^\beta w^{-1}  x = \underline \lambda
\end{equation*}
for all strictly dominant coweights $\beta$, that is,
for dominant
coweights in the interior of the positive Weyl chamber.
\end{lemma}

\begin{proof}
Let us recall some results from \cite{FZ}.
Let $e_i,h_i,f_i$ for $i=1,\dots,r$ be the standard generators of
the Lie algebra $\gg$ of $G$. Denote by
\mbox{$\alpha_i,\dots,\alpha_r\in \hh^*$} the simple roots of
$\gg$. For $p\in \KK$ define
\begin{equation*}
x_i(p) = \exp (p e_i).
\end{equation*}
Then the following commutation relation holds
\begin{equation*}
\tau^\beta x_{i}(p)=x_i(\tau^{\alpha_i(\beta)}p)\tau^\beta
\end{equation*}

Given a reduced word $\ii=(i_1,\dots,i_k)$ of the longest
element~$w_0$ of the Weyl group, it is shown by Fomin and
Zelevinsky in~\cite{FZ} that every element $n$ of $N$ can be
written as a product
\begin{equation}
\label{eq:product} n=x_{i_1}(p_1)\cdots x_{i_k}(p_k)
\end{equation}
with $p_i\in\CC$.
Moreover, they express the $p_i$'s in terms of generalized minors
of $n$.
Since generalized minors can be defined for the group $G(\KK)$ the
same way they were defined for $G$ in \cite{FZ},
it is clear that (\ref{eq:product}) holds for
$n\in N(\KK)$ and~$p_i\in \KK$.

Assume $x=wnw^{-1}\cdot \underline \lambda$, and $n\in N(\KK)$
decomposes as in (\ref{eq:product}); then
\begin{eqnarray*}
\lim_{\tau\to 0} w\tau^\beta w^{-1} x &=& \lim_{\tau\to 0}w\
\tau^\beta
\ x_{i_1}(p_1)\cdots x_{i_k}(p_k)\ w^{-1} \underline \lambda\\
&=&\lim_{\tau\to 0}w\ x_{i_1}(\tau^{\alpha_{i_1}(\beta)}p_1)\cdots
x_{i_k}(\tau^{\alpha_{i_k}(\beta)}p_k) \ \tau^\beta
w^{-1}\underline \lambda.
\end{eqnarray*}
If $\beta $ is  a strictly dominant coweight, then all
$\alpha_{i_m}(\beta)$ are positive integers, so that all
$x_{i_m}(\tau^{\alpha_{i_m}(\beta)}p_m)$ approach the identity as
$\tau$ goes to zero. This shows that the above limit is equal
to~$\underline \lambda$.

This finishes the proof, since for a fixed $w$, $\Gr$ is the union
of all $S^w_\lambda$.
\end{proof}

\subsection{Moment map images of strata}
As explained in \cite{AP}, the loop Grassmannian can be thought of
as an infinite-dimensional symplectic manifold with a Hamiltonian
action of the maximal compact torus $T_K \subset T$. In
particular, it is possible to define a moment map $\momm$ on $\Gr$
with range in the Lie algebra, $\lie(T_K)$. (Although the usual
codomain for moment maps is the dual $\lie(T_K)^*$, here it is
most naturally $\lie(T_K)$, the vector space in which the
coweights of $G$ lie.  This is because moment map images are
closely related to the representation theory of the Langlands dual
group~\cite{MV,VG}, whose weights are identified with the
coweights of $G$.) We do not reproduce the definition of the
moment map, but state the only two properties of this map from
\cite{AP} which we are going to use:
\begin{enumerate}
\item {The moment map image of the fixed
point $\underline \lambda$ of the $T$ action  is the corresponding
coweight~$\lambda$.}
\item {\label{prop2}For a dominant coweight $\lambda$, the moment
map image of both $G_\lambda$ and $\overline {B_\lambda}$ is the
convex hull of the Weyl group orbit containing $\lambda$.}
\end{enumerate}
Note that every piece of each of the decompositions of $\Gr$ is
$T$-invariant; indeed, this will be true of every variety we
consider.

One of our main interests is in the MV-polytopes, which are moment
map images of MV-cycles.  To compute these, we will need the
moment map images of the $S^w_\lambda$, which we now describe.
Let~$C_\lambda$ be the cone inside~$\lie(T_K)$ with
vertex~$\lambda$ and spanned by the negatives of the simple roots;
note that $\alpha\leq \lambda$ if and only if $\alpha\in
C_\lambda$. Let~$C_\lambda^w$ be the cone produced by acting on
$C_\lambda$ by the Weyl group element $w$.

\begin{lemma}
\label{lem:cone}
The moment map image of $S_\lambda^w$ is $C^w_\lambda$.
\end{lemma}

\begin{proof}

It suffices to prove the statement for $S_\lambda$.
Using property (\ref{prop2}) of the moment map and the fact
that~\mbox{$\sigma_{-\alpha}B_{\lambda+\alpha}\subseteq S_\lambda$}
for any dominant coweight $\alpha$, we see that $\momm(S_\lambda)$
contains $C_\lambda$.
Conversely, suppose $x\in S_\lambda$.  Then the closure
$\overline{T\cdot x}$ of the torus orbit through $x$
is contained in $\overline{S_\lambda}$.
By \cite{Ati82},
$\momm(\overline{T\cdot x})$
is a convex polytope.
Every vertex $\beta$ of this polytope is a coweight contained
in $\momm(\overline{S_\lambda})$;
hence $\beta\leq \lambda$ and so $\beta \in C_\lambda$.
So $C_\lambda$ contains the polytope, which contains
$\momm(x)$.
\end{proof}

\subsection{Moment map images of algebraic subvarieties}

By an algebraic subvariety of $\Gr$, we mean any algebraic
subvariety of one of the $\overline {G_\lambda}$. Our interest in
the following theorem is that it applies to MV-cycles.

\begin{theorem}
\label{thm:cones} The moment map image of a compact irreducible
torus-invariant algebraic subvariety~$Z$ of~$\Gr$ is a polytope
given by intersecting cones $C_{\lambda_w}^w$, one cone for every
element of the Weyl group.
\end{theorem}

\begin{proof}
The moment map image of every such variety is a convex polytope~\cite{Bri}.

Fix $w$ in the Weyl group. The intersection $Z\cap S^w_\lambda$,
if nonempty, contains some torus orbit; since~$Z$ is closed, it
follows by Lemma~\ref{lem:limit} that $\underline\lambda\in Z$. So
there are only finitely many coweights $\lambda_1,\dots,\lambda_k$
for which this intersection is nonempty, and all are contained in
$\momm(Z)$; therefore $\momm(Z) \subseteq \bigcup_i
C^w_{\lambda_i}$ by Lemma~\ref{lem:cone}. Among these coweights,
there must exist one, $\lambda_w$, that is greater than all the
others; for otherwise there would be two of them such that the
line segment joining them was not contained in $\momm(Z)$,
contradicting its convexity.

We have just shown that $\momm(Z)$ lies inside the
convex polytope $P=\bigcap_w S^w_{\lambda_w}$ and contains the
$\lambda_w$'s.
To prove equality of the polytopes,
it remains only to show that $P$
has no vertices besides the~$\lambda_w$'s.
Let $\beta$ be any vertex of $P$
and fix a hyperplane through $\beta$ intersecting $P$ only in $v$
and generic in the sense that it contains no root direction.
This hyperplane determines a set of positive roots and
a positive root cone $C$.  Near $\beta$, $P$ is the
intersection of some set of half-spaces,
each on one side of a hyperplane spanned by roots;
therefore $P$ must be contained in $C$.  Since $C$ is a
translation of $C^w_{\lambda_w}$ for some $w$,
we must have $\beta=\lambda_w$.
\end{proof}

\begin{remark}
As explained in Section~\ref{sec:lattice-model}, every connected
component of $\lgrgl$ is isomorphic to $\lgrsl$. Hence
Theorem~\ref{thm:cones} also holds for $\gln$, even though it is
not simply-connected or semisimple.
\end{remark}

\section{Moment map images of torus orbits}
\label{sec:moment-images}

In this section we study moment map images of torus orbits---in
particular of orbits through lattices strongly compatible to a
given Kostant picture. This allows us to compute the moment map
images of the cycles $\overline{\wcomp(\pp,\lambda)}$ and prove
Theorems~\ref{thm:main} and~\ref{thm:main2}.

\subsection{Moment maps of general torus orbits}

For a lattice $Y$, denote by $P(Y)$ the moment map image of the
closure of the torus orbit through~$Y$. We now describe the
vertices of $P(Y)$.

Recall that for a
subset $I$ of $\{1,\dots,n\}$, we defined $Y_I=Y\cap V_I$ as the
intersection of $Y$ with the columns $I$. Further define
$d_I(Y)=\dimz(Y_I)$.
Given any permutation $w$ of $\{1,\dots n\}$,
define the coweight $\mu^w(Y)=(\mu_1,\dots,\mu_n)$ by
\begin{equation*}
\mu_{w(i)}=\base_{w(i)}(Y)+d_{w([i..n])}(Y)-d_{w([i+1..n])}(Y),
\end{equation*}
where $[i..n]$ is short for $\{i,i+1,\dots,n\}$.  Equivalently,
$\mu_{w(i)}$ is the largest $j$ such that $t^{-j}e_{w(i)}$ is in
$\proj_{V_{w(i)}}\big(Y\cap V_{w([i..n])}\big)$.

\begin{lemma}
\label{lem:moment}$P(Y)$ is a convex polytope whose vertices
are~$\mu^w(Y)$.
\end{lemma}

\proof By Theorem~\ref{thm:cones} and Lemmas~\ref{lem:limit}
and~\ref{lem:cone}, it is enough to show that
\begin{equation*}
\lim_{\tau\to 0} w\tau^\beta w^{-1} Y = \underline{\mu^w(Y)}
\end{equation*}
for any strictly dominant coweight
$\beta=(\beta_1>\dots>\beta_n)$, where the action of $\tau^\beta$
on a lattice is given by multiplying each $V_i$ by
$\tau^{\beta_i}$. Using Weyl group symmetry, it suffices to prove
it for the identity permutation $\id$ only.

By definition of $\mu^\id(Y)$ there exists a basis of $Y$ such
that exactly $k_i=d_{[i..n]}(Y)- d_{[i+1..n]}(Y)$ vectors
$y^i_1,\dots, y^i_{k_i}$ of this basis lie in $V_{[i..n]}$ and
have linearly independent projections onto~$V_i$. Moreover, since
$\beta=(\beta_1>\dots>\beta_n)$, we have
\begin{equation*}
\lim_{\tau\to 0} \tau^\beta  \langle y^i_m\rangle  = \lim_{\tau\to
0} \tau^\beta  \langle \tau^{-\beta_i}y^i_m\rangle = \langle
\proj_{V_{i}}y^i_m\rangle.
\end{equation*}
Therefore
\begin{equation*}
\lim_{\tau\to 0} \tau^\beta  Y = Y_0 \oplus \lim_{\tau\to 0}
\tau^\beta \langle y^i_m\rangle_{1\leq i\leq n, 1\leq m\leq k_i} =
Y_0\oplus\langle \proj_{V_{i}}y^i_m\rangle_{1\leq i\leq n, 1\leq
m\leq k_i}=\underline{\mu^\id(Y)}. \ \ \ \ \ \square
\end{equation*}

\begin{remarks}

\noindent (1) Specifying the numbers $\base_i(Y)$ and $d_I(Y)$ for
all $I$ is equivalent to specifying the polytope $P(Y)$. Indeed,
Lemma~\ref{lem:moment} gives the vertices in terms of them.
Conversely, these numbers are determined by the vertices:
$\base_j(Y)$ equals $\mu_{w(n)}$ for any choice of $w$ with
$w(n)=j$; then
$d_I(Y)=\sum_{i=k}^n \big(\mu_{w(i)} - \base_{w(i)}(Y)\big)$ for
any choice of $w$ with $w([k..n])=I$. The very interesting
question of which
polytopes arise is equivalent to the question of which
combinations of numbers $d_I(Y)$ are possible.

\smallskip

\noindent (2) Specifying the numbers $\base_i(Y)$ and $d_I(Y)$ for
any interval $I=[\ell..r]$ is by definition equivalent to
specifying the pair $(\highweight(Y),\kst(Y))$.

\smallskip

\noindent (3) If $P$ is a polytope, let $M_P$ denote the set of
lattices $Y$ for which $P=P(Y)$. This gives the decomposition of
$\lgr$ into torus orbit types, and it is now easy to see that it
is a refinement of the decomposition into the $\wcomp(\pp,\lambda)$'s.
Indeed, if $Y$ is a lattice and we set $P=P(Y)$, $\pp=\kst(Y)$,
$\lambda=\highweight(Y)$ then we conclude by remarks (1) and (2)
that $M_P \subseteq \wcomp(\pp,\lambda)$.

\smallskip

\noindent (4) For each element $w$ of the Weyl group $W$ there is a
decomposition $\wcomp^w(\pp,\lambda)$ analogous to
$\wcomp(\pp,\lambda)$, where $\wcomp(\pp,\lambda)=
\wcomp^{\rm{id}}(\pp,\lambda)$: for a lattice $Y$ let $w\cdot Y$
be the lattice obtained by reordering the columns according to $w$
and set $\wcomp^w(\pp,\lambda) = \{ Y \mid \pp=\pp(w\cdot Y)
\mbox{ and } \lambda=\mu^w(Y) \}$. Then the intersection of these
$n!$ decompositions is the decomposition into torus orbit
types. Indeed, since every subset $I$ of $\{1,\dots,n\}$ is mapped
onto an interval $[\ell..r]$ by some Weyl group element, we see by
remarks (1) and (2) that
$\wcomp_{P(Y)} = \bigcap_{w\in W}
\wcomp^w(\pp(w\cdot Y), \mu^w(Y)).
$
\end{remarks}

\subsection{Moment images of torus orbits through weakly
compatible lattices}

Given a lattice~$Y$, let us call $\mu^\id(Y)$ and $\mu^{w_0}(Y)$
the \emph{highest} and \emph{lowest} vertices of $P(Y)$.

\begin{lemma}
\label{lem:weak-moment}
All polytopes $P(Y)$ for lattices $Y\in \wcomp(\pp,\lambda)$ have
the same highest and lowest vertices.  Moreover the highest vertex
is $\lambda$.
\end{lemma}

\begin{proof} This is an immediate consequence of Lemma~\ref{lem:moment}
and the definitions of $\kst(Y)$ and $\highweight(Y)$.
\end{proof}

\begin{remark}
The lowest vertex is the coweight with $i^{\rm th}$ component
$\lambda_i-l_i+r_i$, where $l_i$ is the number of loops of $\pp$
with left end $i$ and $r_i$ is the number of loops of $\pp$ with
right end $i$.
\end{remark}

\subsection {Moment images of torus orbits through strongly
compatible lattices}

\begin{theorem}
\label{thm:strong-moment} A lattice $Y$ is in $\scomp(\pp,\lambda)$ if
and only if the polytope $P(Y)$ is equal to $P(\pp,\lambda)$.
\end{theorem}

\begin{proof}
Suppose $Y \in \scomp(\pp,\lambda)$. To show $P(Y)=P(\pp,\lambda)$,
we show that these polytopes have the same vertices; that is, for
each permutation $w$, the coweight $\mu^w(Y)=(\mu_1,\dots,\mu_n)$
equals the coweight $\nu(w)=(\nu_1,\dots,\nu_n)$ defined in
Section~\ref{algorithm}. Since the numbers
$\mu_{w(2)},\dots,\mu_{w(n)}$ depend only on the lattice $\hat
Y_{w(1)}$ strongly compatible to $\hat\pp_{w(1)}$, it suffices by
induction to show that $\mu_{w(1)}=\nu_{w(1)}$. Recall that
$\mu_{w(1)} = \base_{w(1)}(Y)+\dimz(Y)-\dimz(\hat Y_i)$ and
$\nu_{w(1)}= \lambda_{w(1)}-l_{w(1)}+N_{w(1)}$; we will check that
the right sides are equal to each other. By definition of
$\lambda=\highweight(Y)$ we have
$\base_{w(1)}=\lambda_{w(1)}-l_{w(1)}$. Since $Y$ is compatible to
$\pp$, we know that $\hat Y_i$ is weakly compatible to $\hat
\pp_i$ by Lemma~\ref{lem:induction}; therefore
$\dimz(Y)-\dimz(\hat Y_i)=|\pp|-|\hat \pp_i|=N_{w(1)}$.

To show the converse, recall that the polytope $P(Y)$ uniquely
determines the numbers $d_I(Y)$.
It follows, by what we just proved, that if
$Y \in \scomp(\pp,\lambda)$ then $d_{\bar I}(Y)=|\hat \pp_\II|$
for any ordered set $\II$ and the corresponding unordered
complement $\bar I$. So we will show that if a lattice $Y\in
\wcomp(\pp,\lambda)$ is not strongly compatible to $\pp$ then
\renewcommand{\labelitemi}{($\star$)}
\begin{itemize} \item
there exists an ordered set $\II$ with $d_{\bar I}(Y)\neq
|\hat\pp_\II|$.
\end{itemize}

Since $Y$ is not strongly compatible to $\pp$, there exists an
ordered set $\JJ$ such that $\hat Y_{J}$ is weakly compatible to
$\hat\pp_\JJ$ but not compatible to it. Obviously, if we can prove
($\star$) for~$\hat Y_J$ instead of for $Y$ it will also follow for
$Y$, so without loss of generality we assume that $Y$ is not
compatible to $\pp$. In particular, there exists a loop $L$
passing through columns $\ell,\dots,r$ for
which condition~(\ref{property'}$'$) on page~\pageref{property'}
fails. Consider the ordered set
$$
\JJ=(1,\dots,\ell-1,n,\dots,r+1).
$$
If $L'$ is the largest loop in $\pp$ encircling the same dots as $L$, then
it is easy to see that $\hat\pp_\JJ=\pp^{L'}$ and $\hat
Y_J=Y_{L'}$. So, if we can prove $(\star)$ for $Y_{L'}$ and
$\pp^{L'}$ instead of for $Y$ and $\pp$, it will also follow for $\pp$.

So, without loss of generality, we can assume that $Y$ is in
$\wcomp(\pp,\lambda)$ but is not compatible to $\pp$ and
that~(\ref{property'}$'$) fails for a column $i$ and a loop $L$
encircling all dots. We claim that $\dimz(\hat Y_i)\neq
|\hat\pp_i|$, which implies $(\star)$ for $\II=(i)$. Indeed, $
\dimz(\hat Y_i)=|\pp|-\dim (\proj_{V_i}(Y)/\proj_{V_i}(Y_0))$ and
$|\hat\pp_i|=|\pp|-N_i$ where~$N_i$ is the number of levels into
which the picture $\pp$ is broken during the collapse along
column~$i$. Notice that loop $L$ is removed during this collapse,
since it is the only loop in its level.
Because~(\ref{property'}$'$) fails for loop~$L$ at column~$i$ we
conclude that $\dim (\proj_{V_i}(Y)/\proj_{V_i}(Y_0))< N_i$; hence
$\dimz(\hat Y_i)>|\hat\pp_i|$ and we have proved $(\star)$.
\end{proof}

\subsection{Proofs of the main theorems.}
\label{sec:proof of the main theorem}

Given two coweights $\alpha$ and $\beta$, recall that MV-cycles
are defined to be the irreducible components of the closure of the
set  $S_{\alpha,\beta}=S_\alpha\cap S^{w_0}_\beta$.
By Lemmas~\ref{lem:limit}~and~\ref{lem:cone}, $S_{\alpha,\beta}$
contains all lattices $Y$ for which the highest and lowest
vertices of $P(Y)$ are $\alpha$ and $\beta$ respectively.
Therefore, by Lemma~\ref{lem:weak-moment}, $S_{\alpha,\beta}$ is a
union of several pieces of the decomposition $\wcomp(\pp,\alpha)$.
By Proposition~\ref{prop:M-is-smooth}, all these pieces are
Zariski open sets of the same dimension, $\len(\pp)$; hence each
closure $\overline{\wcomp(\pp,\alpha)}$  is an MV-cycle. On the
other hand, for any $\wcomp(\pp,\lambda)$, there exists a $\beta$
for which $\wcomp(\pp,\lambda)\subseteq S_{\lambda,\beta}$, by
Lemma~\ref{lem:weak-moment}. This proves Theorem~\ref{thm:main}.
By Propositions~\ref{prop:compatible}
and~\ref{strongly-compatible} the closure of $\scomp(\pp,\lambda)$
is the same as the closure of $\wcomp(\pp,\lambda)$. By
Theorem~\ref{thm:strong-moment} the moment map image of
$\overline{\scomp(\pp,\lambda)}$ and hence of
$\overline{\wcomp(\pp,\lambda)}$ is $P(\pp,\lambda)$. This proves
Theorem~\ref{thm:main2}.

\end{document}